\documentclass[12pt, reqno]{amsart}
\usepackage{amssymb, amsthm, amsmath, amsfonts, bbm}
\usepackage{enumerate}
\usepackage{array}	
\usepackage{lineno}
\usepackage{color}

\usepackage{hyperref}
\usepackage{esvect}

\usepackage{anysize}
\usepackage[pdftex]{graphicx}

\usepackage[makeroom]{cancel} %crossing out a formula

\newcommand*\patchAmsMathEnvironmentForLineno[1]{%
  \expandafter\let\csname old#1\expandafter\endcsname\csname #1\endcsname
  \expandafter\let\csname oldend#1\expandafter\endcsname\csname end#1\endcsname
  \renewenvironment{#1}%
     {\linenomath\csname old#1\endcsname}%
     {\csname oldend#1\endcsname\endlinenomath}}% 
\newcommand*\patchBothAmsMathEnvironmentsForLineno[1]{%
  \patchAmsMathEnvironmentForLineno{#1}%
  \patchAmsMathEnvironmentForLineno{#1*}}%
\AtBeginDocument{%
\patchBothAmsMathEnvironmentsForLineno{equation}%
\patchBothAmsMathEnvironmentsForLineno{align}%
\patchBothAmsMathEnvironmentsForLineno{flalign}%
\patchBothAmsMathEnvironmentsForLineno{alignat}%
\patchBothAmsMathEnvironmentsForLineno{gather}%
\patchBothAmsMathEnvironmentsForLineno{multline}%
}

  \evensidemargin
\oddsidemargin
\newcommand{\argmax}{\mathop{\mathrm{argmax}}\nolimits}
\newcommand{\argmin}{\mathop{\mathrm{argmin}}\nolimits}

\newcommand{\Var}{\mathop{\mathrm{Var}}\nolimits}

\newcommand{\supp}{\mathop{\mathrm{supp}}\nolimits}

\newcommand{\sgn}{\mathop{\mathrm{sgn}}\nolimits}

\makeatletter
\DeclareFontFamily{U}{tipa}{}
\DeclareFontShape{U}{tipa}{m}{n}{<->tipa10}{}
\newcommand{\arc@char}{{\usefont{U}{tipa}{m}{n}\symbol{62}}}%

\newcommand{\arc}[1]{\mathpalette\arc@arc{#1}}

\newcommand{\arc@arc}[2]{%
  \sbox0{$\m@th#1#2$}%
  \vbox{
    \hbox{\resizebox{\wd0}{\height}{\arc@char}}
    \nointerlineskip
    \box0
  }%
}
\makeatother

\newcommand{\ubar}[1]{\underline{#1}}

\DeclareMathOperator{\per}{per}
\DeclareMathOperator{\length}{length}
\DeclareMathOperator{\area}{area}

\DeclareMathOperator{\conv}{conv} 
 
\DeclareMathOperator{\im}{im}

\DeclareMathOperator{\cl}{cl}
\DeclareMathOperator{\diam}{diam}

\DeclareMathOperator{\rint}{rint}
\DeclareMathOperator{\intr}{int}
\DeclareMathOperator{\epi}{epi}

\DeclareMathOperator{\id}{Id}

\newtheorem{cor}{Corollary}[section]
\newtheorem{prop}{Proposition}[section]
\newtheorem{lem}{Lemma}[section]
\newtheorem{theorem}{Theorem}[section]
\newtheorem*{theorem*}{Theorem}

\theoremstyle{remark}
\newtheorem{remark}{Remark}[section]

\newtheorem{example}{Example}[section]

\def\R{\mathbb{R}}

\def\N{\mathbb{N}}
\def\E{\mathbb{E}}
\def\P{\mathbb{P}}

\def\S{\mathbb{S}}
\def\I{\mathbbm{1}}

\newcommand{\rank}{\mathop{\mathrm{rank}}\nolimits}

\newcounter{fig}
\newcommand{\f}{\refstepcounter{fig} Fig. \arabic{fig}. }

\begin{document}

\ifpdf
\DeclareGraphicsExtensions{.pdf, .jpg, .tif, .mps}
\else
\DeclareGraphicsExtensions{.eps, .jpg, .mps}
\fi

\title[Large deviations of convex hulls of planar random walks]{Large deviations of convex hulls \\of planar random walks and Brownian motions}

\author{Arseniy~Akopyan}
\address{Arseniy V.\ Akopyan, Institute for Information Transmission Problems RAS} %, Bolshoy Karetny per.\ 19, Moscow,  127994, Russia}
\email{akopjan@gmail.com}

\author{Vladislav~Vysotsky}
\address{Vladislav Vysotsky, University of Sussex and St.\ Petersburg Department of Steklov Mathematical Institute}
\email{v.vysotskiy@sussex.ac.uk}

\thanks{This paper was written when AA was affiliated to IST Austria. His work was supported by the People Programme (Marie Curie Actions) of the European Union's Seventh Framework Programme (FP7/2007-2013) under REA grant agreement n$^\circ$[291734] and European Research Council (ERC) under the European Union's Horizon 2020 research and innovation programme (grant agreement No 78818 Alpha). This paper was partially written when VV was affiliated to Imperial College London, where his work was supported by the People Programme (Marie Curie Actions) of  the European Union's Seventh Framework Programme (FP7/2007-2013) under REA grant agreement n$^\circ$[628803]. 
%VV is supported in part by the RFBR Grant 19-01-00356.
}

\begin{abstract}
We prove large deviations principles (LDPs) for the perimeter and the area of the convex hull of a planar random walk with finite Laplace transform of its increments. 

We give explicit upper and lower bounds for the rate function of the perimeter in terms of the rate function of the increments. These bounds coincide and thus give the rate function for a wide class of distributions which  includes the Gaussians and the rotationally invariant ones. For random walks with such increments, large deviations of the perimeter are attained by the trajectories that asymptotically align into line segments. However, line segments may not be optimal in general. 
%Here the rate function is convex, but we also give examples with non-convex rate functions. 

Furthermore, we find explicitly the rate function of the area of the convex hull for random walks with rotationally invariant distribution of increments. 
For such walks, which necessarily have zero mean, large deviations of the area are attained by the trajectories that asymptotically align into half-circles. For random walks with non-zero mean increments, we find the rate function of the area for Gaussian walks with drift. Here the optimal limit shapes are elliptic arcs if the covariance matrix of increments is non-degenerate and parabolic arcs if otherwise.

The above results on convex hulls of Gaussian random walks remain valid for convex hulls of planar Brownian motions of all possible parameters. Moreover, we extend the LDPs for the perimeter and the area of convex hulls to general L\'evy processes with finite Laplace transform.

%Finally, we prove counterparts to the above results for L\'evy processes that %have finite Laplace transform, with the rate functions available explicitly for  %Brownian motions of all possible parameters.
\end{abstract}


\maketitle

%\tableofcontents

\section{Introduction} \label{Sec: intro}
%\subsection{Motivation} 
Let $(S_k)_{k \ge 1}$, where $S_k = X_1 + \ldots + X_k$, be a planar random walk with independent identically distributed increments $X_1, X_2, \ldots.$ We assume that the expectation of $X_1$ exists and is finite, and put $\mu := \E X_1$. We are interested in the perimeter $P_n$ and the area $A_n$ of the convex hull $C_n:=\conv(0, S_1, \dots, S_n)$ of the first $n$ steps of the random walk, including the origin. Here, by definition, the perimeter of a line segment  is its doubled length.  
%It is intentional that we do not exclude one-dimensional random walks (whose increments are supported on an affine line in the plane) from our consideration. 

All of our results remain valid for the convex hulls $\conv(S_1, \dots, S_n)$ but it is more natural to consider hulls of the form $C_n$, which allow remarkably simple formulas for their expected perimeters and areas. In fact, Spitzer and Widom~\cite{SpitzerWidom} proved\footnote{\cite{SpitzerWidom} proved formula \eqref{eq: Spitzer-Widom} under the assumption $\P(u \cdot S_k=0)=0$ for every $k \in \N$ and non-zero $u \in \R^2$,  which can be dropped using a simple approximation argument based on the fact that the perimeter is a continuous functional on the space of compact convex sets equipped the Hausdorff distance. This observation also applies to equality \eqref{eq: E area}.} that
\begin{equation} \label{eq: Spitzer-Widom}
\E P_n = 2 \sum_{k=1}^n \frac{\E |S_k|}{k},
\end{equation}
where by $|\cdot|$ we denote the Euclidean norm. This implies that $\E P_n/n \to 2 |\mu|$ as $n \to \infty$, by the law of large numbers and uniform integrability of $(S_k/k)_{k \ge 1}$. 
%(the latter property holds by uniform integrability of $(\frac1k \sum_{i=1}^k |X_i|)_{k \ge 1}$). 
Moreover, $P_n/n \to 2|\mu|$ a.s.\ by McRedmond and Wade~\cite{McRedmondWade}. 

Wade and Xu~\cite{WadeXu} showed (developing the ideas introduced by Snyder and Steele~\cite{SnyderSteele}) that if $\mu \neq 0$ and $\E|X_1|^2<\infty$, then $\Var(P_n)/n \to 4 \sigma_\mu^2/|\mu|^2$, where  $\sigma_\mu^2:=\E(\mu \cdot ( X_1 - \mu ))^2$ and `$\cdot$' denotes the scalar product. Here $\sigma_\mu^2>0$ unless 
%$X_1$ is supported on the line $\{u \in \R^2: \mu \cdot (u - \mu) =0\}$, in %which case 
the trajectory of $(S_k)_{k \ge 1}$ is the graph of a zero-mean one-dimensional random walk. With the exception of this degenerate case,  the variance of the perimeter grows linearly (when $\E|X_1|^2<\infty$), and moreover, the sequence $(P_n)_{n \ge 1}$ satisfies a central limit theorem for $\mu \neq 0$ (see~\cite{WadeXu}) and a limit theorem under the scaling $n^{-1/2}$  for $\mu = 0$ (see Wade and Xu~\cite{WadeXu2}). The latter result  follows naturally from the invariance principle using the continuous mapping theorem. The degenerate case $\mu \cdot (X_1 - \mu) =0$ a.s.\ is more tricky. Alsmeyer et al.~\cite{AKMV} proved that $\Var(P_n) =O(\log n)$ under $\E |X_1|^3<\infty$; it is likely that $\Var(P_n)$ may grow super-logarithmically (contradicting the corresponding conjecture in~\cite{WadeXu}) when $\E |X_1|^3=\infty$ and $\E |X_1|^2<\infty$, but a proof is still missing. Yet there is no central limit theorem for $(P_n)_{n \ge 1}$, although Alsmeyer et al.~\cite{AKMV} established the ones for the lengths of the convex minorant and the concave majorant of $n$-step one-dimensional random walks (the sum of these quantities is $P_n$; their variances may grow polynomially when $\E |X_1|^3=\infty$ but not much known about their correlation).

The Spitzer--Widom formula \eqref{eq: Spitzer-Widom} admits various generalizations to higher dimensions, including explicit formulas for the expected mean width, surface area, volume, and other intrinsic volumes of the convex hulls, see Barndorff-Nielsen and Baxter~\cite{BNielsen} or Vysotsky and Zaporozhets~\cite{VysotskyZaporozhets}. In particular, for the area of the convex hull of a planar random walk, 
\begin{equation} \label{eq: E area}
\E A_n = \frac12 \sum_{\substack{j, k \ge 1 \\ j + k \le n }} \frac{\E|\det[ S_j, S_k' ] |}{jk},
\end{equation}
where $(S_k')_{k \ge 1}$ is an independent copy of $(S_k)_{k \ge 1}$. Furthermore, the invariance principle naturally implies (see Wade and Xu~\cite{WadeXu2}) that if $\E |X_1^2| <\infty$, then the sequence $(A_n)_{n \ge 1}$ satisfies a limit theorem under the scaling $n^{-3/2}$ for $\mu \neq 0$ and $n^{-1}$ for $\mu=0$. 

%Let $\Sigma:=\E (X_1-\mu) (X_1- \mu)^\top$ be the covariance matrix of %$X_1$ and put $\sigma^2:= \tr \Sigma$; then $\sigma^2= \E |X_1-\mu|%^2$. Assume that $\sigma^2 < \infty$. With a certain effort, equality~%\eqref{eq: E area} implies, by the central limit theorem, that $\E A_n/ %n^{3/2} \to \frac13 |\mu| \sqrt{2 \pi(\sigma^2 - \sigma_\mu^2)}$ for $%\mu \neq 0$ and $\E A_n /n \to \frac{\pi}{2} \sqrt{\det \Sigma}$ for $%\mu=0$,  as $n \to \infty$ (cf.~Wade and Xu~\cite{WadeXu2}). 

In this paper we study large deviations probabilities for the perimeter and the area of the convex hull of the random walk. This describes very atypical behaviour of these quantities, as opposed to the results above on their typical behaviour. In particular, we will consider the logarithmic asymptotics of $\P(P_n \ge 2 x n)$ for $x> |\mu|$ and $\P(P_n \le 2 x n)$ for $x<|\mu|$, and $\P(A_n \ge a n^2)$ for $a > 0$. We will also describe the limit shape of the trajectories, scaled by the factor of $n^{-1}$ in both time and space (this explains the scalings of $P_n$ and $A_n$),  resulting in such large deviations probabilities.

To the best of our knowledge, there is only one rigorous result in this direction. Snyder and Steele~\cite{SnyderSteele} obtained the following non-sharp concentration inequality for the perimeter for random walks with bounded increments: if $|X_1| \le M$ a.s.\ for some $M >0$, then 
\[
\P(|P_n - \E P_n| \ge xn) \le 2\exp(- x^2 n/(8 \pi^2 M^2)), \quad x \ge 0.
\]
Claussen et al.~\cite{Claussen+} gave a numerical analysis of atypically large values of the perimeter and the area of the convex hull and concluded that these quantities ``seem ...\ to obey a large deviations principle'' ({\it LDP}, in short) for random walks with standard Gaussian increments. 
%However, the heuristic explanation of this conclusion given in~\cite{Claussen+} is not convincing. 
There are few follow-up numerical papers on related questions by the same group of authors. 

The other results on atypical behaviour of convex hulls include the works by Khoshnevisan~\cite{Khoshnevisan} and Kuelbs and Ledoux~\cite{KuelbsLedoux}, who considered a.s.\ superior limits of monotone functionals of convex hulls (including the perimeter and the area) of zero-mean finite-variance random walks and  standard Brownian motions scaled as in the law of iterated logarithm.

%considering large deviations of the following quantities: the perimeter and the area of the joint convex hulls of several %planar random walks by Dewenter et al.~\cite{Dewenter+}, the volume and the surface area of the convex hull of walks in %higher dimensions, and even those 

% followed by the work of Dewenter et al.~\cite{Dewenter+} which concerns joint convex hulls of several random walks, 

The simulation-based conclusions of~\cite{Claussen+}  in fact easily follow (see Section~\ref{Sec: main proofs} below) from the contraction principle applied to Mogulskii's LDP for trajectories of random walks with finite Laplace transform of their increments. The main task is to obtain {\it explicitly} the rate functions in these LDP's in terms of the rate function of the increments. For the perimeter, we found the rate function (Corollary~\ref{cor: thm holds}) for a wide class of random walks (see Proposition~\ref{prop: conjecture holds}) including all Gaussian walks (this is not an expected result at all), and also gave the upper and the lower bound  valid for general walks with finite Laplace transform of  increments (Theorem~\ref{thm: perimeter LD}). For the area, we found the rate function for random walks that have rotationally invariant distributions of increments with finite Laplace transform (Theorem~\ref{thm: LD area}) and for Gaussian random walks with arbitrary drift (Theorem~\ref{thm: LD Gaussian area} and Proposition~\ref{prop: area 1D}).  In all these results we identified the asymptotic form of optimal trajectories of the walk resulting in the large deviations.

Furthermore, we extended the above results on random walks, which have increments in discrete time, to convex hulls of planar L\'evy processes (Theorem~\ref{thm: LD Levy}) with finite Laplace transform, including Brownian motions. Convex hulls of general L\'evy processes were studied e.g.\ by Molchanov and Wespi~\cite{MolchanovWespi}.

Lastly, we extended the LDPs for $P_n$ and $A_n$ to random walks whose increments have Laplace transform finite only in a neighbourhood of zero  (Proposition~\ref{prop: Cramer}). 

\medskip
The paper is organized as follows. In Section \ref{sec: notation} we introduce notation and in particular, define the \emph{radial minimum rate function} and state its properties. In Sections~\ref{sec: results perimeter} and \ref{sec: results area} we present our main results on large deviations for the perimeter and the area of the convex hull of a planar random walk. The continous-time counterparts are given in Section~\ref{sec: continuous time}. Further generalizations are discussed in Section~\ref{sec: further}. 
In Section~\ref{Sec Convex} we prove basic properties of the radial minimum rate function for general increments of the walk and its convexity for Gaussian walks. %in Section~\ref{Sec: convexity ubar_I}. 
Section \ref{Sec: basics on LDP} contains essentials on large deviations relevant to this paper.
%, then in \ref{Sec: main proofs} we prove that it holds true for perimeters and %area of a convex hull.
The proofs of our LDPs for the perimeter and the area (for random walks), including  computations of the rate functions, are given in Section~\ref{Sec: main proofs}. The core parts of these computations are unified by the use of geometric inequalities of isoperimetric type. The proofs for L\'evy processes are in Section~\ref{sec: continuous time proof}. Finally, in Section~\ref{sec: LDP Cramer} we give a partial result for walks with Laplace transform of increments finite only in a neighbourhood of zero.

\section{Main results}

\subsection{Notation} \label{sec: notation}
Recall that the {\it Legendre--Fenchel transform} or the {\it convex conjugate} of a function $F:\R^d \to \R \cup \{+ \infty\}$ (where $d \ge 1$) with a non-empty {\it effective domain} $\mathcal{D}_F:=\{u: F(u) < \infty\}$  is the function  $F^*:\R^d \to \R \cup \{+ \infty\}$ defined by
\[
F^*(v):= \sup_{u \in \R^d} \bigl( u \cdot v - F(v) \bigr), \quad v \in \R^d.
\] 
The conjugate function $F^*$ is convex and lower semi-continuous on $\R^d$; $F$ itself does not need to be convex.  
%By convention, we always mean that a convex function defined of a convex subset of $\R^d$ is extended to  the whole space and equal %infinity outside their effective domains. 
%By saying that a function on $\R^d$ is convex we always mean that its domain is $\R^d$ (which is not necessarily its effective %domain). 
Recall that any convex function $F$ is continuous on the relative interior $\rint \mathcal{D}_F$ of its effective domain (Rockafellar~\cite[Theorem~10.1]{Rockafellar}) so the property of lower semi-continuity is needed to characterize $F$ only near the relative boundary of $\mathcal{D}_F$. By $\conv F$ we denote the {\it largest convex minorant} or the {\it convex hull} of $F$, i.e.\ the convex function with the epigraph $\conv(\epi F)$, which is a subset of $\R^{d+1}$. Thus, we use the  notation ``$\conv$'' both for functions and sets.

The {\it cumulant generating function} $K(u):=\log \E e^{u \cdot X_1}$ is convex by Jensen's inequality and satisfies $K(0)=0$. %and even strictly convex unless $X_1$ belongs to a line a.s.,  
Its convex conjugate $I:=K^*$ is the {\it rate function} of $X_1$. 
%$$I(v):=K^*(v)= \sup \limits_{u \in \R^2} \bigl(u \cdot v - \log \mathcal{L} (u) \bigr), \quad v \in \R^2,$$ 
This function satisfies $I(\mu)=0$ and is non-negative, lower semi-continuous, and continuous on $\rint \mathcal{D}_I$, where $\rint$ stands for the {\it relative interior} (taken in the induced topology of the affine hull of  $\mathcal{D}_I$). In the main results of this paper (namely, the LDPs for $P_n$ and $A_n$) we assume that the Laplace transform of the increments $\mathcal{L}(u):= \E e^{u \cdot X_1}$ is finite for all $u\in \R^2$. For example, this is trivially true when the support of $X_1$ is bounded. Under this assumption, $K$ is infinitely differentiable on $\R^2$ and $I$ is {\it strictly} convex on its effective domain $\mathcal{D}_I$; see Barndorf-Nielsen \cite[Corollary~7.1]{B-N} and Vysotsky~\cite[Corollary~1]{VysotskyStrict}. 

The effective domain of $I$ is known to satisfy (see~\cite[Proposition~1]{VysotskyStrict})
\begin{equation} \label{eq: D_I inclusions}
\rint(\conv( \supp(X_1)))  \subset \mathcal{D}_{I} \subset  \cl(\conv( \supp(X_1))),
\end{equation}
where $\supp(X_1)$ is the topological support of the distribution of $X_1$. 
%A complete description of $\mathcal D_I$ is available in~\cite{VysotskyStrict}. 
Furthermore, put
\[
r_{min}:= \inf \{|u|: u \in \conv( \supp(X_1)) \}, \quad r_{max}:= \sup \{|u|: u \in  \supp(X_1) \}.
\]
Note that $r_{min} \le |\mu| \le r_{max}$, where the second inequality is strict unless $X_1 = \mu$ a.s.\ and the first inequality is strict unless $\mu \cdot (X_1 - \mu) = 0$ a.s. 

On occasions, we will give general statements assuming that the random vector $X_1$ takes values in $\R^d$ with an arbitrary $d \ge 1$ rather than merely in $\R^2$. With no  risk of confusion, in such cases in $\mathcal{L}(u)= \E e^{u \cdot X_1}$ we take $u \in \R^d$  and understand $I$, $K$, etc.\ accordingly. Then we will usually assume that $X_1$ satisfies merely the {\it Cram\'er moment assumption} $0 \in \intr \mathcal D_{\mathcal L}$.

Define the {\it radial maximum} and {\it radial minimum} functions
\begin{equation} \label{eq: radial functions}
\ubar{I}(r):= \inf_{\ell \in \S^{d-1}} I(r\ell), \quad \bar K(p):= \sup_{\ell \in \S^{d-1}} K(p\ell),  \quad p, r \ge 0,
\end{equation}
where $\S^{d-1}$ stands for the unit sphere in $\R^d$ centred at $0$, and put $\ubar I(r):= \infty$ and $\bar K(p) := \infty$ for $p, r<0$. 
%We will further use this upper and lower bar notation for arbitrary functions on $\R^d$. 
Note that the function $\ubar I$ admits the following geometric interpretation: the epigraph of $\ubar I(|v|)$ is the union of all rotations of the epigraph of $I(v)$ about the vertical axis. Clearly, the supremum and the infimum above are always attained at some points since the Laplace transform is continuous, $I$ is lower semi-continuous, and spheres in $\R^d$ are compact. Thus, the respective sets of minimal and maximal directions
\begin{equation} \label{eq: extreme directions}
\ubar{\Lambda}_r:=\argmin \limits_{\ell \in \S^{d-1}} I(r \ell) , \quad \bar{\Lambda}_p:= \argmax \limits_{\ell \in \S^{d-1}}  K(p \ell), \quad r,p \ge 0
\end{equation}
are always non-empty. Note that the $\argmax$ will not change if we replace $K$ by $\mathcal L$ in \eqref{eq: extreme directions}.

%Careful with Prop. 1d)!
%For convenience of notation, we put $\bar{\Lambda}_0(\mathcal{L}):=\mu/|\mu|$ if $\mu\neq %0$ and $\bar{\Lambda}_0:=\mu/|\mu|$.  
%Without any risk of confusion, throughout Section~\ref{Sec: intro} we will use the short notation $\ubar{\Lambda}_r$ and $%\bar{\Lambda}_p$ for $\ubar{\Lambda}_r(I)$ and $\bar{\Lambda}_p(\mathcal{L})=\bar{\Lambda}_p(K)$, respectively. 

Recall that a point $u$ in a convex set $C \subset \R^d$ is called {\it extreme} if there is no way to express $u = \alpha u_1 + (1- \alpha) u_2$ for some $u_1, u_2 \in C$ and $\alpha \in (0,1)$ except by taking $u_1=u_2=u$. Every extreme point of $C$ belongs to the {\it relative boundary} $\partial_{\text{rel}} C$ of $C$, defined by  $\partial_{\text{rel}} C: =C \setminus \rint C$. An extreme point $u$ of a convex set $C \subset \R^d$ is called {\it exposed} if $C \cap L = \{ u \}$ for some hyperplane $L$ supporting~$C$.  

The radial minimum rate function $\ubar{I}$ and the sets of minimal directions $\ubar{\Lambda}_r$ appear in most of our results on the perimeter of the convex hull. Let us state some of their properties. Let us agree that by $[|\mu|, r_{max}]$ we will mean the half-line $[|\mu|, \infty)$ if $r_{max}=\infty$. 

\begin{lem} \label{lem: properties of I_}
Assume that $X_1$ is a random vector in $\R^d$, where $d \ge 1$, such that $0 \in \intr \mathcal D_{\mathcal L}$. 
 %$\E e^{u \cdot X_1}< \infty$ for any $u\in \R^d$.
\begin{enumerate}[a)]
\item \label{item: I_ eff domain} The effective domain $\mathcal{D}_{\ubar{I}}$ of $\ubar I$ is an interval that satisfies $\intr \mathcal{D}_{\ubar{I}} = (r_{min}, r_{max})$;
\item \label{item: I_ unimodal} The function $\ubar I$ is lower semi-continuous; satisfies $\ubar I (|\mu|)=0$; is strictly decreasing and convex (also strictly if $\mathcal{D}_{\mathcal L }=\R^d$) on $[r_{min}, |\mu|]$; and is strictly increasing on $[|\mu|, r_{max}]$;
%\item \label{item: I_ eff domain} If $r_{max} < \infty$, then $\mathcal{D}_{\ubar{I}}$ is an interval such that $\partial \mathcal{D}%_{\ubar{I}} =\{r_{min}, r_{max}\}$, otherwise $\mathcal{D}_{\ubar{I}}$ is either $(r_{min}, \infty)$ or $[r_{min}, \infty)$;
%\item If $r \in \partial \mathcal{D}_{\ubar{I}}$ and $\ubar I(r)<\infty$, then $\ubar{I}(r) = - \log %\P(X_1 = r \ell)$ for any $\ell \in \ubar{\Lambda}_r$ and $\P(|X_1| > r) \cdot \P(|X_1| < r) = 0$.
%either $\P(|X_1| > r) = 0$ or $\P(|X_1| < r) = 0$ (unless $X_1 = \mu$ a.s., when both equations hold true);
\item \label{item: I_ discontinuous} Suppose that $\ubar I$ is discontinuous at a point $x \in [|\mu|, r_{max}]$. Then for any $\ell \in \ubar{\Lambda}_x$, $x \ell$ is an exposed point of $\mathcal D_I$ and $\ubar I(x) = -\log \P(X_1 = x \ell)<\infty$.
\item \label{item: one element} For any $ r \in (r_{min}, |\mu|]$, the set $\ubar{\Lambda}_r$ contains a unique element, which we denote by $\ell_r$.
\end{enumerate}
\end{lem}
%In addition, we give explicit formulas for $\ubar I (r_{min})$ and $\ubar I (r_{max})$ presented below in Lemma~\ref{lem: boundary %values}.

Combining Part~\ref{item: I_ discontinuous} with the second inclusion in \eqref{eq: D_I inclusions} gives: 
\begin{cor} \label{cor: I_ continuous}
$\ubar I$ is continuous on $(r_{min},\infty)$ if $\P(X_1=u)=0$ for any $u \in \partial_{\text{rel}}(\conv(\supp(X_1)))$.
\end{cor}

We stress that the function $\ubar I$ may be discontinuous on its effective domain, and may be non-convex even if it is continuous; see 
Remark~\ref{rem: counter example} in Section~\ref{sec: results perimeter}  and Example~\ref{ex: I_ discontinuous}  in Section~\ref{Sec: Basic Proof}.

%Let $H$ be a compact subset of $C([0,1], \R^2)$, the space of continuous $\R^2$-valued functions equipped with the %topology of uniform convergence. We say $H$ is a {\it set of asymptotic trajectories} of the random walk $(S_k)_{k \ge 1}$ %on a sequence of events $E_k \in \sigma(S_1, \ldots, S_k)$, $k \ge 1$, of positive probability if for any $\varepsilon %>0$, 
%\begin{equation} \label{eq: def optimal}
%\lim_{n \to \infty} \P \Bigl(\max_{1 \le k \le n} \Bigl | \frac{S_k}{n} - h\Big(\frac{k}{n} \Big) \Bigr| \le \varepsilon %\text{ for some } h \in H \Bigl | E_n \Bigr) = 1.
%\end{equation}
%We will use this definition in concrete cases; in general, such a set may not exit. We call $H$ {\it minimal} if %\eqref{eq: def optimal} ceases to hold for $H$ replaced by each of its closed proper subsets.

\subsection{Large deviations of the perimeter} \label{sec: results perimeter}
Denote by $AC_0([0,1];\R^2)$ the set of coordinate-wise absolutely continuous functions $h$ on $[0,1]$ such that $h(0)=0$. We will occasionally refer to functions from $[0,1]$ to $\R^2$ as (planar) {\it curves} or {\it trajectories}. Denote by $\im(\cdot)$ the image of a function, that is the  set of its values as the argument varies over the effective domain. Let $P(C)$ denote the perimeter of a non-empty convex set $C \subset \R^2$, so $P_n=P(C_n)$. 

We now state our first main result.

%Writing $I$ in the polar coordinates, for any $r \le |\mu|$ it holds that
%$$\ubar I(r)= \inf_{0 \le \theta \le 2 \pi} I(r, \theta) = \min_{(\rho, \theta) \in \mathcal{D}_I: \atop \rho \le r, 0 \le \theta \le %2 \pi} I(\rho, \theta),$$ hence $\ubar I$ is convex on $[0, |\mu|]$ as the minimum of a convex function over a convex set.

\begin{theorem} \label{thm: perimeter LD}
Assume that $X_1$ is a random vector in $\R^2$ such that $\mathcal D_{\mathcal L} = \R^2$.

1.\ The sequence $(P_n/(2n))_{n \ge 1}$ satisfies the LDP in $\R$ with speed $n$ and the tight rate function 
\begin{equation} \label{eq: rate func I_P}
\mathcal J_P(x):=\min_{\substack{h \in AC_0([0,1]; \R^2): \\ P(\conv(\im h)) = 2x }} \int_0^1 I(h'(t)) dt.
\end{equation}
This function shares the properties of $\ubar I$ stated in Parts \ref{item: I_ eff domain}  and \ref{item: I_ unimodal}  of Lemma~\ref{lem: properties of I_}.

2.\ We have $\mathcal J_P = \ubar I$ on $[0, |\mu||]$ and $\conv \ubar I \leqslant \mathcal J_P \leqslant \ubar I$ on $[|\mu|, \infty)$, and a bit more: %Moreover, 
\begin{equation} \label{eq: LD per < E}
\lim_{n \to \infty} \frac 1n \log \P(P_n \le 2 xn) = 
\begin{cases}
-\ubar I(x), & x \in ( r_{min}, |\mu|]\\
\log \P(|X_1|=r_{min}), & x = r_{min}
\end{cases}
\end{equation}
and for any $x \in [|\mu|, r_{max}]$, 
\begin{equation} \label{eq: LD per > E}
-\ubar I (x) \le \liminf _{n \to \infty} \frac 1n \log \P(P_n \ge 2 xn)  \le \limsup _{n \to \infty} \frac 1n \log \P(P_n \ge 2 xn)  \le -\conv \ubar I (x).
\end{equation}

3.\ For any $\varepsilon > 0$, we have
\begin{equation}
\label{eq: LD shape <} 
\lim_{n \to \infty} \P \Bigl(\max_{0 \le k \le n} \Bigl | \frac{S_k}{n} - \frac{k}{n} x \ell_x \Bigr| \le \varepsilon \Bigl | \Bigr. P_n \le 2 xn \Bigr) = 1, \qquad x \in ( r_{min}, |\mu|],
\end{equation}
where $\ell_x$ was defined in Lemma~\ref{lem: properties of I_}.\ref{item: one element}, 
\begin{equation}
\label{eq: LD shape >} 
\lim_{n \to \infty} \P \Bigl(\max_{0 \le k \le n} \Bigl | \frac{S_k}{n} - h(k/n) \Bigr|  \le \varepsilon \text{ for some } h \in H_P(x) \Bigl | \Bigr. P_n \ge 2 xn \Bigr) = 1, \qquad x \in [|\mu|,  r_{max}),
\end{equation}
where $H_P(x)$ denotes the set of minimizers in \eqref{eq: rate func I_P}. If $\ubar I(x) = \conv \ubar I (x)$ for an $x \in (r_{min}, r_{max})$, then $H_P(x)=\{t \mapsto t x \ell\}_{\ell \in \ubar \Lambda_x}$.
\end{theorem}

%The lower bound in \eqref{eq: LD per > E} is slightly stronger than the lower bound in the LDP in Part~1 combined with the %inequality $J_P \leqslant \ubar I$; this is typical for one-dimensional LDPs. 

We refer to the elements of the sets $H_P(x)$ as the {\it optimal trajectories} (for the perimeter).

Let us give a few comments. The limit shape results \eqref{eq: LD shape <} and \eqref{eq: LD shape >} mean that if $\ubar I(x) = \conv \ubar I(x)$, then large deviations of the perimeter are attained on trajectories that asymptotically align into line segments and move with constant speed. Note that under  $\ubar I(x) = \conv \ubar I(x)$, equality \eqref{eq: LD shape >} does not assert that every direction in $\ubar{\Lambda}_x$ can be attained.
%; this claim is not accessible using the rough logarithmic asymptotics given above. 
This equality means that the epigraph of $\ubar I$ admits a support line at the point $(x, \ubar I(x))$. This is true for every $x$ iff $\ubar I$ is convex, in which case the rate function in the LDP for the perimeter is $\mathcal J_P = \ubar I$. In Proposition~\ref{prop: properties of I_} below we will provide a tractable condition, stated directly in terms of the Laplace transform of increments, for checking the equality $\ubar I(x)= \conv \ubar I(x)$ for a given $x$. Moreover, this proposition relates the sets of optimal directions $\ubar{\Lambda}_x$ in \eqref{eq: LD shape >} to more tractable sets $\bar \Lambda_p$ defined in terms of the Laplace transform; cf.~\eqref{eq: extreme directions}.

The idea of our proof of Theorem~\ref{thm: perimeter LD} is as follows. Equalities \eqref{eq: rate func I_P}, \eqref{eq: LD shape <}, and \eqref{eq: LD shape >} follow from an LDP for trajectories of random walks combined with the contraction principle. If $\ubar I(x) = \conv \ubar I(x)$, an additional geometric argument yields that the set $H_P(x)$ consists of curves of minimal length with the fixed perimeter $2x$ of their convex hull. A known geometric result (Corollary~\ref{cor:length of curve} in the Appendix) asserts that the image of such a curve is a line segment of length $x$. This curves must move with constant speed by strict convexity of $I$.

With this geometric optimality property of line segments, it is tempting to assume that $\mathcal J_P = \ubar I$  and \eqref{eq: LD shape >} always hold true. However, in general, 
\begin{center}
{\it the optimal trajectories are not necessarily linear.} 
\end{center}
Hence it may be that $\mathcal J_P \neq \ubar I$, as shown in Example~\ref{ex: counter}, which follows the next remark.

%This is demonstrated in Example~\ref{ex: non-straight} of Section~\ref{sec: %non-straight}, which is constructed using the fact that there are distributions  %(with $\mathcal D_{\mathcal L} \neq \R^2$) with non-strictly convex rate %functions.

\begin{remark} \label{rem: counter example}
For $\mu \neq 0$, the upper bound in Part 2 of Theorem~\ref{thm: perimeter LD} can improved  to $\mathcal J_P \leqslant \ubar {I_0}$ on $[|\mu|, \infty)$, where $I_0:=\conv(\I_{\R \setminus \{0\}} \cdot I)$ is the function with the epigraph $\conv(\epi I \cup \{0\})$. This follows from considering the set of trajectories 
\[
D :=\big \{ h \in AC_0[0,1]: h'=v \text{ on } [0,s], h'=\mu \text{ on } [s,1] \text{ for some } s \in (0,1), v \in T \big \},
\]
where $T:=\{ v \in \mathcal D_I: I(tv) \ge t I(v) \text{ for any } t \in \R\}$, whose energies satisfy 
\[
\int_0^1 I(h')dt = s I(v) + (1-s) I(\mu) = sI(v) < I(sv).
\]
Note in passing that it is easy to check that $T=(K \circ \nabla I)^{-1}(0) $ if $\nabla I$ is defined on $\R^2$. 
%since $sI(v)-I(sv)$ is a concave function of $s$ that attains its maximum value %zero at $s=1$ by 
%$$
%\frac{\partial}{\partial s} \big( sI(v)-I(sv)\big) \Big|_{s=1} \Big.= I(v) - \nabla %I(v)  \cdot v  = - K(\nabla I(v)) = 0.
%$$

We actually have $\mathcal J_P = \ubar {I_0}$ on $[|\mu|, \infty)$ if the distribution of $X_1$ is supported on the straight line $\mu \R$ (and satisfies $\mathcal D_{\mathcal L} = \R^2$). Even in this degenerate case $\ubar{ I_0}$ differs from $\ubar I$ on $[|\mu|, \infty)$ if there is an $a>1$ such that $K(-a\mu)=0$ and $I(a\mu) > I(-a\mu)$ (cf.~Example~\ref{ex: counter}). There are only two directions, therefore both optimal trajectories, one of which belongs to $D$,  start {\it moving backwards} at some moment. This is very counter-intuitive! In this case $\mathcal J_P$ is non-convex, and so is $\ubar I$. Unfortunately, the case $\P(X_1 \in \mu \R)=1$ is the only type of distribution of the increments  with a possibly non-convex $\mathcal J_P$ where we found $\mathcal J_P$ explicitly.

		\begin{center}
		\includegraphics{graph-1.mps}		
		
		\f \label{fig: counterexample} $\ubar{I_0} (x) < \ubar I(x)$ for $x \in [12.77,21.88]$.
	\end{center}	

\begin{example}[{\it Non-linear optimal trajectories}] \label{ex: counter}
Consider the distribution of $X_1$ given by the mixture of the uniform distribution on the line segment $[11,13] \times \{0\}$ taken with weight $49/50$ and the delta distribution at $(-38,0)$ taken with weight $1/50$; see Figure~\ref{fig: counterexample}.
%, where the blue, orange, and green graphs represent, respectively, $I(x,0)$, $I(-x,0)$, and the tangent line to  $I(-x,0)%$ passing through zero, for $x \ge 0$.  %Here $T=\{-0.093,11\}$.
\end{example}
\end{remark}

We now turn our attention to the function $\conv \ubar I$. The following result gives a simple description of $\conv \ubar I$ directly in terms of $\bar K$ and a condition when it equals $\ubar I$ at a given point.
%$\bar K$, the radial maximum of cumulant generating function. 
%The use of this statement is illustrated by the fact that the only direct proof of the convexity of $\ubar I$ that we have  for %Gaussian distributions is very cumbersome even in this basic case.
%The use of this statement is illustrated by the fact that a direct proof of the convexity of $\ubar I$ is very cumbersome even for %Gaussian distributions. 
Denote by $(\cdot)'_+$ and $(\cdot)'_-$ respectively right and left derivatives of a function of real argument.

\begin{prop} \label{prop: properties of I_}
Assume that $X_1$ is a random vector in $\R^d$, $d \ge 1$, such that $\mathcal D_{\mathcal L} = \R^d$. Then 
\begin{enumerate}[a)]
\item \label{item: conv I_ =} %\label{item: K^- properties} 
$\bar K$ is an increasing convex function on $[0, \infty)$ satisfying $\bar K'_+(0)=|\mu|$ and 
$$
\conv \ubar I = (\bar{K})^* \text{ on }[|\mu|, \infty);
$$
\item \label{item: I_=conv I_} If $r \in \cl (\im(\bar{K}'))$ (and $r \ge |\mu|$), then $\ubar I(r) = \conv \ubar I(r)<\infty$;
\item \label{item: differentiable} For any $p \in (0,\infty)$, the one-sided derivatives satisfy
$$\bar K_+'(p) = \max_{\ell \in \bar{\Lambda}_p } |\nabla K(p \ell)| \quad \text{and} \quad \bar K_-'(p) = \min_{\ell \in \bar{\Lambda}_p } |\nabla K (p \ell)|.$$
%Importantly, for any $r \ge |\mu|$,
%\item \label{item: conv I_ =} $\conv \ubar I = (\bar{K})^*$ on $[|\mu|, \infty)%$;
%The condition that $\ubar I(r) = \conv \ubar I(r)<\infty$ and $(r, \ubar I(r))$ is an extremal %point of $\epi (\conv \ubar I)$ is equivalent to  $r \in \intr (\im(\bar{K}'))$ (the image is open if we show that $\bar K$ is strictly convex);
\item \label{item: Lambdas} If  $p \in (0,\infty)$ and $r \ge |\mu|$ are such that $\bar{K}'(p)=r$, then $\ubar{\Lambda}_r = \bar{\Lambda}_p$.
\end{enumerate}
\end{prop}

The main result of the proposition is Part~\ref{item: conv I_ =}, which relates the radial maximum of the logarithmic Laplace transform $K$ to the radial minimum of its convex conjugate $I$. This assertion is actually a general fact valid for arbitrary convex functions; see Proposition~\ref{prop: generalization} of Section~\ref{Sec: radial proofs}, which yields a stronger version of Proposition~\ref{prop: properties of I_} under the Cram\'er moment assumption. Part~\ref{item: I_=conv I_} is an easy consequence of Part~\ref{item: conv I_ =} and the well-known fact that the Legendre--Fenchel transform maps kinks of a convex function (in our case, $\bar K$) into linear segments of its convex conjugate. Part~\ref{item: differentiable}, which clarifies the possible reason of non-differentiability of $\bar K$, follows by a standard application of the method of Lagrange multipliers. Part~\ref{item: Lambdas}, which follows naturally from Parts~\ref{item: conv I_ =}-\ref{item: differentiable}, claims that the slowest directions of $I$ are exactly the fastest directions of $K$ (equivalently, of $\mathcal{L}$) at the corresponding radii.

The main use of Proposition~\ref{prop: properties of I_}  is through its following corollaries.

\begin{cor} \label{cor: L differentiable}
$\ubar I$ is strictly convex on $\mathcal D_{\ubar I}$ if $\bar K$ is differentiable on $(0, \infty)$. 
\end{cor}
\begin{cor} \label{cor: condition for diff}
$\bar K$ is differentiable if there exists a continuous mapping $\ell: (0, \infty) \to \S^{d-1}$ such that $\ell(p) \in \bar \Lambda_p$ for any  $p \in (0,\infty)$.  
\end{cor}

These claims follow easily from Parts~\ref{item: I_=conv I_} and~\ref{item: differentiable}, respectively, using convexity of $\bar K$; see~Section~\ref{Sec: radial proofs}. We  do not assume that $\ell(p)$ is differentiable, otherwise  Corollary~\ref{cor: condition for diff} becomes trivial. 

We now present a few types of distributions with convex radial minimum rate function~$\ubar I$. By $\| \cdot \|$ we denote the the largest eigenvalue of a symmetric real matrix.
%Denote by $\lambda_1$ the largest eigenvalue of $\Sigma$.
%Denote by $\|\Sigma\|:=\sup_{\ell: |\ell|=1} |\Sigma \ell|$ the operator norm of $\Sigma$ induced by the Euclidean norm. 

\begin{prop} 
\label{prop: conjecture holds}
Assume that $X_1$ is a random vector in $\R^d$, $d \ge 1$. Then the function $\ubar I$ is convex in either of the following cases:
\begin{enumerate}[a)]
\item \label{item: pseudo-affine} 
$X_1=A Y_1 + \mu$, where $Y_1$ is a random vector in $\R^k$, $1 \le k \le d$, with rotationally invariant distribution and $A$ is a $d \times k$ real matrix such that $A A^\top \mu =\|A A^\top \! \| \mu$;
%$\max_{u \in |\mu| \S^{d-1}} |A^\top u|=|A^\top \mu|$;
%$A A^\top \mu = \lambda_1 \mu$ for $\lambda_1:=\sup_{u \in \S^{d-1}} |A^\top u|$.

%$X_1= \Sigma^{1/2} Y_1 + \mu$, where $Y_1$ is a random vector in $\R^d%$ with rotationally invariant distribution satisfying $0 \in \intr \mathcal %D_{\mathcal L_{Y_1}}$, and it holds that $\mu^\top \Sigma \mu = \lambda_1 |%\mu|^2  $;
\item \label{item: Gaussian} $X_1$ is Gaussian. %$(\mu, \Sigma)$.
\end{enumerate}
\end{prop}
\begin{cor} \label{cor: thm holds}
For the above types of distributions, we have $\mathcal J_P = \ubar I$ (cf.~Theorem~\ref{thm: perimeter LD}).
\end{cor}

We  regard that Corollary~\ref{cor: thm holds} gives the rate function  $\mathcal J_P $ explicitly, since finding the radial minimum $\ubar I$ is a standard optimization problem (which is much simpler than~\eqref{eq: rate func I_P}) solvable  using the method of Lagrange multiplies.

%For $\mu \neq 0$, the assumption on the matrix $A$ is satisfied iff $\mu$ %is a maximal eigenvector of  $A A^\top$. 

We now comment on Proposition~\ref{prop: conjecture holds}. Note that in Case~\ref{item: pseudo-affine},  the matrix $A A^\top$ is proportional to the {\it covariance matrix} of $X_1$ given by $\Sigma:=\E (X_1 X_1^\top) - \mu \mu^\top$. The assumption on $A$ is always satisfied if $\mu=0$.

For Case~\ref{item: pseudo-affine}, convexity of $\ubar I$ follows rather directly from that of $I$. Our proof for Case~\ref{item: Gaussian} rests on Corollaries~\ref{cor: L differentiable} and~\ref{cor: condition for diff} and uses properties of quadric curves to construct a path $\ell(p)$. 
%The same approach can be applied for Part~\ref{item: pseudo-affine} but we %prefer to prove it directly using convexity of $I$.

The condition in Corollary~\ref{cor: condition for diff} is trivially satisfied if the set $\cap_{p > 0} \bar{\Lambda}_p$ is non-empty, i.e.\ there exists a direction that maximizes the Laplace transform at all radii. This rather restrictive assumption naturally holds true for {\it either} linearly transformed {\it or} shifted rotationally invariant distributions of increments. Both cases are covered by Case~\ref{item: pseudo-affine} of Proposition~\ref{prop: conjecture holds}, where $\cap_{p > 0} \bar{\Lambda}_p$ is the set of maximal eigenvectors of $\Sigma$ of unit length if $\mu=0$ and  $\cap_{p > 0} \bar{\Lambda}_p = \{\mu /|\mu|\}$  if $\mu \neq 0$. For general {\it affine} transforms (i.e.\ compositions of linear transforms and translations) of rotationally invariant distributions, we were able to prove convexity of $\ubar I$  only for Gaussian distributions, as per Case~\ref{item: Gaussian}. We will see that here $\cap_{p > 0} \bar{\Lambda}_p $ is empty unless the Gaussian distribution of $X_1$ is degenerate or $X_1$ satisfies the assumptions of Case~\ref{item: pseudo-affine}.

%Note that Case~\ref{item: linear} describes all possible linear transforms of $Y_1$: by the singular %value decomposition theorem, for any $2 \times 2$ matrix $A$ it holds that $A Y_1 \stackrel{d}{=} %(AA^\top)^{1/2} Y_1$.

%In fact, Condition~\ref{condition: unique} implies the differentiability by Proposition~\ref{prop: properties of I_}.\ref{item: %differentiable}. For the second one, if there exists a direction $\ell \in \cap_{p > 0} \bar{\Lambda}_p$, which is maximal for all %radii, then $\bar K(p) = \log \mathcal{L}(p \ell), p > 0,$ is differentiable since so is $\mathcal{L}$. 

%Note that Condition~\ref{condition: max dir} is very restrictive. 
%$$\mathcal{L}(u) =1+\mu \cdot u + \frac12 (\E (X_1 X_1^\top)  u, u) + o(|u|^2), \qquad u \to 0,%$$ we have
%Indeed, it is easy to check that
%\begin{equation} 
%\label{eq: Taylor K}
%K(u) = u^\top \mu + \frac12 u^\top \Sigma u + o(|u|^2), \qquad u \to 0,
%\end{equation} 
%hence Condition~\ref{condition: max dir} is equivalent to the assumption that for $\mu \neq 0$, each $\bar \Lambda_p$ contains $\mu / %|\mu| $ and for $\mu =0$, the set $\cap_{p > 0} \bar{\Lambda}_p$ contains a maximal eigenvector of the covariance matrix $\Sigma$ of %$X_1$. 

\subsection{Large deviations of the area}
	\label{sec: results area}
%Let $A_n$ denote {\it area} of the convex hull $C_n$. There is a simple expression for $%\E A_n$ similar to the Spitzer--Widom formula~\eqref{eq: Spitzer-Widom}, see Barndorff%-Nielsen and Baxter~\cite{BNielsen} or Vysotsky and Zaporozhets~%\cite{VysotskyZaporozhets}. Further, the invariance principle naturally implies that $A_n$ %satisfy a limit theorem (under the scaling $n^{-1}$ for $\mu=0$ and $n^{-3/2}$ for $\mu %\neq 0$), as proved by Wade and Xu~\cite{WadeXu2}. Numerical studies of large %deviations probabilities of $A_n$ were performed by Claussen et al.~\cite{Claussen+}, who %arrived to the same conclusion as for the perimeter, namely that the area follows a large %deviations principle.

The main geometric argument used in our proof of Theorem~\ref{thm: perimeter LD} to find the rate function in the LDP for the perimeter can be applied directly to obtain the rate function in the LDP for the area $A_n$ of the convex hull. The problem reduces to finding a planar curve of the unit length that maximizes the area of its convex hull. This question is known as one of Ulam's problems. Although it is very similar to the classical Dido problem and of course has the same answer that the curve is a half-circle  (Moran~\cite{Moran}), it appears that this Ulam problem does not allow an easy solution by reduction. 
%(despite of the opinion of some authors, including Croft et al~\cite[Problem A28]{Croft+}). 
The corresponding isoperimetric inequality  easily yields the following result. 

Denote by $A(C)$ the area of a non-empty convex set $C \subset \R^2$, so $A_n=A(C_n)$. 

\begin{theorem}
\label{thm: LD area}
Assume that $X_1$ is a random vector in $\R^2$ such that $\mathcal D_{\mathcal L} = \R^2$.

1. The sequence $(A_n/n^2)_{n \ge 1}$ satisfies  an LDP with speed $n$ and the tight rate function 
\begin{equation} \label{eq: rate func I_A}
\mathcal J_A(a):=\min_{\substack{h \in AC_0([0,1]; \R^2): \\ A(\conv(\im h)) = a }} \int_0^1 I(h'(t)) dt,
\end{equation}
which is strictly increasing on $\mathcal D_{\mathcal J_A}$ and satisfies $\mathcal J_A(0)=0$. In particular, for any continuity point $a \ge 0$ of $\mathcal J_A$, we have
\begin{equation} \label{eq: LDP A}
\lim_{n \to \infty} \frac 1n \log \P(A_n \ge an^2) = -\mathcal J_A(a).
\end{equation}
The set of minimizers in \eqref{eq: rate func I_A}, denoted by $H_A(a)$, is such that for any $\varepsilon > 0$, 
\begin{equation}
\label{eq: A shape} 
\lim_{n \to \infty} \P \Bigl(\max_{0 \le k \le n} \Bigl | \frac{S_k}{n} - h(k/n) \Bigr|  \le \varepsilon \text{ for some } h \in H_A(a) \Bigl | \Bigr. A_n  \ge an^2 \Bigr) = 1, \quad a \in \intr(\mathcal{D}_{\mathcal J_A}).
\end{equation}

2. Suppose that the distribution of $X_1$ is rotationally invariant. Then $\mathcal J_A(a)=\ubar I (\sqrt{2 \pi a })$; equality \eqref{eq: LDP A} is valid for every $a \ge 0$; and for any $a \in [0, r^2_{max}/(2 \pi))$, 
\begin{equation}
\label{eq: A shape invar} 
H_A(a)= \bigg \{ \sqrt{\frac{2a}{\pi }} \bigl(\cos (\pm \pi t + \alpha) - \cos \alpha, \, \sin (\pm \pi t  + \alpha) - \sin \alpha \bigr) \bigg \}_{\alpha \in \R}.
\end{equation}
\end{theorem}

Thus, for rotationally invariant distributions, large deviations of the area are attained on the trajectories that asymptotically align into half-circles and move with constant speed. Note that for such distributions, $\ubar I$ is convex. 

We will refer to the elements of the sets $H_A(a)$ as the {\it optimal trajectories} (for the area). 

\begin{remark} 
\label{Rem: elliptic}
Assume that the covariance matrix $\Sigma$ of $X_1$ is non-degenerate, satisfies $\mathcal D_{\mathcal L}=\R^2$, and $\Sigma^{-1/2} X_1$ has a rotationally invariant distribution, whose rate function we denote by $I_1$. 
Then $(\Sigma^{-1/2} S_k)_{k \ge 1}$ is a random walk with a rotationally invariant distribution of increments, and the area of its convex hull satisfies
\[
\area(\conv(\Sigma^{-1/2} S_1, \dots, \Sigma^{-1/2} S_n)) = (\det \Sigma)^{-1/2} A_n.
\]
Hence by Theorem~\ref{thm: LD area} we have $\mathcal J_A(a)=\ubar {I_1} \big(\sqrt{2 \pi a/ (\det \Sigma)^{1/2} }\big)$. To rewrite $\mathcal J_A$ in terms of $I$, note that $\ubar {I_1} (|u|) = I_1(u) = I(\Sigma^{1/2} u)$ for $u \in \R^2$, where the last equality follows by changing variables in the definition of $I$. This gives  $\ubar {I_1} (r)= \ubar I (r \sqrt{\lambda_1})$ for $r \ge 0$, where  $0< \lambda_2 \le \lambda_1$  are the eigenvalues of $\Sigma$, 
%Since the set $\ubar \Lambda_x$ contains the unit eigenvector of $\Sigma$ %corresponding to  $\lambda_1$ (which is the unique element of $\ubar\Lambda_x$  if $\lambda_2 < \lambda_1$), 
hence  $\mathcal J_A(a)=\ubar I \big(\sqrt{2 \pi a } (\lambda_1 / \lambda_2)^{1/4}\big).$
%$\mathcal J_A(a)=\ubar I (\sqrt{2 \pi a /\det \Sigma})$ - bylo tak. Eto oshibka???

Moreover, equality \eqref{eq: A shape invar}  remains valid for $a\in [0, \area(\supp(X_1))/(2 \pi^2))$ if we multiply the factor $\sqrt{2a/\pi}$ in \eqref{eq: A shape invar}  by $(\det \Sigma)^{-1/2} \Sigma^{1/2}$. The optimal trajectories are halves of the ellipse $\supp(X_1)$ divided by the lines passing through its centre.  
\end{remark}

For random walks with a shifted (i.e.\ $\mu \neq 0$) rotationally invariant distribution of increments, the limit shapes for the area are not universal, unlike the case of the perimeter. We were able to solve only the Gaussian case. We apply the same approach as in the proof of Theorem~\ref{thm: LD area}. In fact, since Gaussian rate functions are quadratic, computation of the rate function for the area reduces to finding a planar curve of fixed length and fixed endpoints that maximizes the area of its convex hull. Pach~\cite{pach1978isoperimetric} proved that such a curve is a circular arc, as in the Dido problem with fixed endpoints. The corresponding isoperimetric inequality yields the following LDP. 

Let us denote by $u^\bot$  a vector $u \in \R^2$ rotated  $\pi/2$ counterclockwise about the origin. 

\begin{theorem} 
\label{thm: LD Gaussian area}
Suppose that $X_1$ has a shifted standard Gaussian$(\mu, \id)$ distribution on $\R^2$ with a non-zero mean $\mu$. Then 
\begin{equation} \label{eq: J_A Gaussian shifted}
\mathcal J_A(a)=4 a \varphi - \frac12 |\mu|^2 \tan^2 \varphi, \qquad a \ge 0,
%4a \varphi \Bigl (\frac{\varphi \sin^2 \varphi}{2 \varphi - \sin 2\varphi} - 1\Bigr),
\end{equation}
%where $R= \frac{|\mu|}{2 \varphi \cos \varphi}$ and 
where $\varphi \in [0, \pi/2)$ is the unique solution to 
%where $\varphi=H^{-1}(a/|\mu|^2)$ and $H:[0, \pi/2) \to [0, \infty)$ is an increasing bijection given by
\[
%H(\phi):=\frac{2\phi - \sin 2 \phi}{8 \phi^2 \cos^2 \phi}.
\frac{2\varphi - \sin 2 \varphi}{8 \varphi^2 \cos^2 \varphi} = \frac{a}{|\mu|^2};
\]
equality \eqref{eq: LDP A} is valid for every $a \ge 0$; and in the basis $\mu, \mu^\bot$, the set of optimal trajectories~is
\begin{equation}
\label{eq: A shape shifted} 
H_A(a)= \Big \{ \frac{1}{2 \varphi \cos \varphi} \bigl( \sin( 2 \varphi t - \varphi) + \sin \varphi , \pm \cos( 2 \varphi t - \varphi)  \mp \cos \varphi \bigr ) \Big\}, \quad a>0.
\end{equation}
\end{theorem}  

Thus, large deviations of the area of the convex hull for random walks with shifted standard Gaussian increments are attained on the trajectories that asymptotically align into either of the two $\mu$-axially symmetric circular arcs of radius $\frac{|\mu|}{2 \varphi \cos \varphi}$ and angle $2 \varphi$ starting at the origin and ending on the $\mu$-axis. The radius is defined so that for either of the two limit curves, the orthogonal projections of  their velocities to the direction of $\mu$ at times $0$ and $1$ both equal $\mu$.  Note that the shifted standard Gaussian rate function $I(v)=\frac12 |v-\mu|^2$ is not constant on the velocity of the optimal trajectories (which move with constant speed), contrasting the results of Theorems~\ref{thm: perimeter LD} and~\ref{thm: LD area}. 

%Notice the scaling relations for the angle and the radius: $\varphi(\mu, a) = \varphi(k %\mu, k^2 a)$ and $|k|R(\mu, a) = R(k \mu, k^2 a)$ for any $k \neq 0$. 

It is easy to show that the asymptotics in Theorem~\ref{thm: LD area} for the Gaussian case appears as the limit case of Theorem~\ref{thm: LD Gaussian area} as $|\mu| \to 0$ with a fixed $a$: since $\varphi \to \pi/2$, the radius tends to $\sqrt{2a/\pi}$ and the right-hand side of \eqref{eq: J_A Gaussian shifted} tends to $\pi a$, which is $\ubar I (\sqrt{2 \pi a })$ for the standard Gaussian distribution.

\begin{remark} \label{Rem: elliptic arcs}
Using the same argument as in Remark~\ref{Rem: elliptic} above, we can easily check that for $X_1$ following any Gaussian distribution with a non-degenerate covariance matrix $\Sigma$ and non-zero drift $\mu$, it holds that
\[
\mathcal J_A(a)=\frac{4a \varphi}{\sqrt{\det \Sigma}}  - \frac12  (\mu^\top \Sigma \mu)  \tan^2 \varphi, \text{ where } \frac{2\varphi - \sin 2 \varphi}{8 \varphi^2 \cos^2 \varphi} = \frac{a}{\sqrt{\det \Sigma} \cdot \mu^\top \Sigma \mu}, \qquad a \ge 0.
\]
With this uniquely defined $\varphi$, the set of optimal trajectories is given by \eqref{eq: A shape shifted} taken in the basis $\mu, (\Sigma^{-1/2} \mu)^\bot$. 
%replace $\mu^\bot$ by $\Sigma^{1/2} (\Sigma^{-1/2} \mu)^\bot$ to update \eqref{eq: A shape shifted}. 
In this general case the optimal limit shapes are elliptic arcs starting at the origin and ending on the $\mu$-axis.
\end{remark}

For completeness of exposition, we consider shifted degenerate Gaussian distributions. Since these arise as the limit case of non-degenerate Gaussian distributions, we can use Remark~\ref{Rem: elliptic arcs} to get the following result, which we present here without a proof.

\begin{prop} \label{prop: area 1D}
Suppose that $X_1 \stackrel{d}{=}(\mu_1, \mu_2 + \sigma Y_1)$, where $\mu=(\mu_1, \mu_2)$, $\mu_1$ and $\sigma$ are non-zero, and $Y_1$ is a standard Gaussian random variable. Then $\mathcal J_A(a) = 6 a^2 \mu_1^{-2} \sigma^{-2}$ for $a \ge 0$;
%\mathcal J_A(a) = -\frac{6 a^2}{|\mu|^2 \tr (\Sigma) - \mu^\top \Sigma %\mu}, \qquad a \ge 0;
\eqref{eq: LDP A} is valid for every $a \ge 0$; and the optimal trajectories are the parabolas
\begin{equation*}
%\label{eq: A shape shifted degen} 
H_A(a)= \big \{\mu t \pm 6a \mu_1^{-1}(0, t- t^2) \big \}, \quad a > 0.
\end{equation*}
\end{prop}

In the case $\mu_1=1$, the proposition describes large deviations for the area of the convex hull of the graph of one-dimensional random walk with Gaussian$(\mu_2, \sigma^2)$ increments. The assumptions of Proposition~\ref{prop: area 1D}  ensure that the distribution of $X_1$ is not supported on the line~$\mu \R$ passing through the origin. 
%Large deviations of the area of the convex hull are attained on two $\mu$-axially symmetric parabolas $\mu t \pm 6a \mu_1^{-1}(0, t- %t^2)$ starting at the origin and ending on the $\mu$-axis. 

\subsection{Convex hulls of L\'evy processes} \label{sec: continuous time}

The above results on convex hulls of random walks, which have increments in discrete time, have the following counterparts in continuous time. 

Assume now that $(S_t)_{t \ge 0}$ is a {\it L\'evy process} on the plane, that is a stochastic process with stationary independent increments and c\`adl\`ag trajectories (i.e.\ right-continuous and having left limits)  taking values in $\R^2$. Then $(S_t)_{t \in \N}$ is a random walk. Conversely,  every random walk with an {\it infinitely divisible} distribution of increments (that for every $n\in \N$ is the $n$-fold convolution of some  distribution) can be regarded as such time-discretization of a L\'evy process $(S_t)_{t \ge 0}$; see Bertoin~\cite[Theorem~I.1]{Bertoin}. 

Consider the convex hull $\mathsf C_T:= \conv(\{S_t\}_{0 \le t \le T})$ for $T > 0$. Its perimeter $\mathrm P_T:=P(\mathsf C_T)$ and area $\mathrm A_T:=A(\mathsf C_T)$ are random variables (see the Appendix). 

The following general result extends Parts 1 of Theorems~\ref{thm: perimeter LD} and~\ref{thm: LD area}.

\begin{theorem} \label{thm: LD Levy}
Assume that $(S_t)_{t \ge 0}$ is a {\it L\'evy process} on the plane such that $S_1=X_1$ and  $\mathcal D_{\mathcal L} = \R^2$. Then the variables
$(\mathrm P_T/(2T))_{T>0}$ and $(\mathrm A_T/T^2)_{T>0}$ satisfy the LDP's in $\R$ (as $T \to \infty$) with speed $T$ and the respective rate functions $\mathcal J_P$ and $\mathcal J_A$ (given in \eqref{eq: rate func I_P} and~\eqref{eq: rate func I_A}).
\end{theorem}

\begin{remark} \label{rem: Bm}
Assume that $(S_t)_{t \ge 0}$ is a planar Brownian motion starting at zero (which is the only L\'evy process with continuous trajectories). By Corollary~\ref{cor: thm holds} we have $\mathcal J_P = \ubar I$, and $\ubar I$ can be found using the method of Lagrange multipliers as in the proof of   Proposition~\ref{prop: conjecture holds}.\ref{item: pseudo-affine} (where we essentially found a similar quantity~$\bar K$).  $\mathcal J_A$ is given explicitly in Theorem~\ref{thm: LD Gaussian area}, Proposition~\ref{prop: area 1D}, and Remarks~\ref{Rem: elliptic},~\ref{Rem: elliptic arcs}, which together cover all possible values of the drift $\mu$ and the covariance $\Sigma$ of~$S_1$.  

Moreover, the limit shape results~\eqref{eq: LD shape <},~\eqref{eq: LD shape >}, \eqref{eq: A shape}  remain valid with $k$ and~$n$ considered as positive reals (instead of integers, as before) and $P_n$, $A_n$ replaced respectively by $\mathsf P_n$, $\mathsf A_n$; this describes the optimal trajectories of the Brownian motion that result  in the large deviations of $\mathsf P_n$ and $\mathsf A_n$.  Indeed, the proofs of Theorems~\ref{thm: perimeter LD} and~\ref{thm: LD area} can be carried over without any changes if instead of Mogulskii's LDP (see Section~\ref{Sec: basics on LDP}) for trajectories  of random walks we apply Schilder's LDP (\cite[Theorem~5.2.3]{DemboZeitouni}) for trajectories  of a Brownian motion in $\R^d$ (this result is stated in~\cite{DemboZeitouni} for a standard Brownian motion but we can convert it into the LDP for $S$ applying the contraction principle to the mapping $f \mapsto \Sigma^{1/2} f(t) + \mu t$ for $f \in AC_0[0,1]$).
\end{remark}

We will prove Theorem~\ref{thm: LD Levy} by reduction to the random walks case, showing that the trajectory of $(S_t)_{0 \le t \le T}$  stays close  to that of $(S_t)_{t \in \N \cup \{0\}}$. This also allows one to extend our simplest limit shape result~\eqref{eq: LD shape <} to general   L\'evy processes. It is also tempting to provide counterparts to~\eqref{eq: LD shape >} and \eqref{eq: A shape}. However, it appears that arguing by reduction to random walks would require additional technical assumptions. Therefore, we do not give any results in this direction. Note that we cannot directly prove counterparts to~\eqref{eq: LD shape <},~\eqref{eq: LD shape >},~\eqref{eq: A shape}, as we did for Brownian motions in Remark~\ref{rem: Bm} above, since  
we are not aware of any LDP for general   L\'evy processes appropriate for the purpose.

\subsection{Further extensions} \label{sec: further}
\subsubsection{Higher dimensions}
One can further consider large deviations of surface area, volume, etc.\ for convex hulls of random walks in higher dimensions. The expected values of these quantities are available through the explicit formulas of~\cite[Section 4]{VysotskyZaporozhets} which generalize~\eqref{eq: Spitzer-Widom}. However, currently we cannot obtain any progress  even for rotationally invariant  distributions of increments. In fact, according to Tilli~\cite{Tilli}, the problem of finding the shape of a curve in $\R^d$, where $d \ge 3$, of unit length that maximizes volume of its convex hull is yet solved only in the class of curves convex in the sense of Schoenberg (i.e., those that intersect no hyperplane at more than $d$ points) and there is no complete solution. Croft et al.~\cite[Problem A28]{Croft+} mention that there are no results on the similar problem of maximizing the surface area, and we are unaware of any progress in this direction.

\begin{remark} \label{rem: thm mean width}
On the other hand, our results for the perimeter of planar random walks can be easily extended for {\it mean width} of convex hulls in higher dimensions, defined in~\eqref{eq: mean width}. A closely related quantity is the {\it first intrinsic volume} of the convex hull, which equals (see~\cite[Eq.\ (14.7)]{SchneiderWeil}) the mean width divided by $\frac{2 v_{d-1}}{d v_d}$, which is mean width of a unit segment,  where $v_d$ denotes volume of a unit ball in $\R^d$.

Assume now that $S_n$ is a random walk in $\R^d$ satisfying $\mathcal D_{\mathcal L} =\R^d$. 
%and extend in the obvious way the definition of $\ubar I$ and the other quantities related %to the Laplace transform $\E e^{u \cdot X_1}$. Lemma~\ref{lem: properties of I_} of %course holds in this setting. 
Denote by $W_n$ and $V_n$ the mean width and first intrinsic volume, respectively, of the convex hull $C_n$. The Spitzer--Widom formula~\eqref{eq: Spitzer-Widom} remains valid (see~\cite[Corollary 3]{VysotskyZaporozhets}) in any dimension if we replace the perimeter $P_n$ of the convex hull $C_n$ by its doubled first intrinsic volume $2 V_n$. Accordingly, our Theorem~\ref{thm: perimeter LD} remains valid if we replace $P_n$ by $\frac{d v_d}{v_{d-1}} W_n = 2 V_n$ so the probabilities change to $\P(W_n \ge \frac{2v_{d-1}}{d v_d} x n)$ and $\P(W_n \le  \frac{2v_{d-1}}{d v_d} x n)$ or, equivalently, more elegant expressions $\P(V_n \ge x n)$ and $\P(V_n \le x n)$.  
The only difference in the proof is that Remark~\ref{rem: mean width} in the Appendix should be used instead of Corollary~\ref{cor:length of curve}.
\end{remark}

\subsubsection{Weaker exponential moments assumptions} The Cram\'er moment assumption $0 \in \intr \mathcal D_{\mathcal L}$ for the increments is a standard minimal requirement to work with large deviations of random walks. However, in the case $\mathcal D_{\mathcal L} \neq \R^2$ we have to regard trajectories of the walk as random elements of the space of functions of bounded variation equipped with a Skorokhod topology  (either of them is strong enough for our purposes). Essentially, this is due to the fact that the rate function $I$ is not super-linear at infinity, and in particular, the infima in~\eqref{eq: rate func I_P} and~\eqref{eq: rate func I_A} may not be attained on absolutely continuous functions. 
%Additional difficulties arise since trajectories of such random walks satisfy a %form of LDP which is weaker than the standard one; see Vysotsky~%\cite{VysotskyNote}. 

The only available large deviations result for such trajectories is the non-standard LDP by Borovkov and Mogulskii~\cite{BorovkovMogulskii2}. It can be applied to our problems  using the contraction principle by Vysotsky~\cite{VysotskyNote}, which yields LDP's for the perimeter and the area in the case $0 \in \intr \mathcal D_{\mathcal L}$; see Proposition~\ref{prop: Cramer} in Section~\ref{sec: LDP Cramer}. The rate functions there are rather complicated but remarkably, they are exactly the same as in the main case $\mathcal D_{\mathcal L} = \R^2$ when $\ubar I$ is convex (for the perimeter) or the distribution of $X_1$ is rotationally invariant (for the area). We do not identify  the optimal trajectories in Proposition~\ref{prop: Cramer}  since there are too many cases to analyse.

\section{Properties of the radial minimum rate function $\ubar I$} \label{Sec Convex}
\subsection{Basic facts from convex analysis} \label{Sec Convex basics} Suppose that $F: \R^d \to \R \cup \{+\infty \}$ is any function with a non-empty effective domain $\mathcal{D}_F$. 

$\bullet$ By \cite[Theorem 12.2 and Corollary 12.1.1]{Rockafellar}) it holds that
\begin{equation} 
\label{eq: involution}
%F^* = (\cl (\conv F))^*, \quad 
F^{**}= \cl (\conv F),
\end{equation} 
where $\cl ( \cdot)$ denotes the {\it closure} of a function, that is the function with the epigraph $\cl (\epi (\cdot))$. Recall that $F$ is lower semi-continuous iff $F = \cl F$, i.e.\ its epigraph $\epi F$ is closed in $(\R \cup \{+\infty\})\times \R^d$ (\cite[Theorem~7.1]{Rockafellar}). Thus the Legendre--Fenchel transform is an involution on the set of lower semi-continuous convex functions. 
%Note that $\cl (\conv F)$ is the largest convex lower semi-continuous function dominated by $F$.

$\bullet$ The {\it Fenchel inequality} $F(u) + F^*(v) \ge u \cdot v$, which holds for any $u, v \in \R^d$, immediately follows from the definition of convex conjugation. 

%If $u$ and $v$ are such that $F(u) + F^*(v) = u \cdot v$, then $\nabla F(u)= v$ if $F$ is %differentiable at $u$ and $\nabla F^*(v) = u$ if $F^*$ is differentiable at $v$. Conversely, if  $%\nabla F(u)= v$, then $F^*(v) = u \cdot v - F(u)$ by the Fenchel inequality since the function $u %\mapsto u \cdot v - F(u)$ is concave. Similarly, $\nabla F^*(v)= u$ implies that $F(u) = v \cdot u - %F^*(v)$. 

$\bullet$ If the function $F$ is convex, then it is continuous on $\rint (\mathcal{D}_F)$ (\cite[Theorem~10.1]{Rockafellar}) and if, in addition, $F$ is lower semi-continuous, then it is continuous on every closed interval contained in $\mathcal D_F$ (\cite[Corollary~7.5.1]{Rockafellar}). If $F$ is convex, finite, and differentiable on an open convex set $C \subset \R^d$, then $F$ continuously differentiable on $C$ (\cite[Corollary 25.5.1]{Rockafellar}). If $C=\R^d$, then $F^*$ is {\it strictly convex} on $\rint (\mathcal{D}_{F^*})$, which means that $F^*$ is linear on no line segment with the endpoints in $\rint (\mathcal{D}_{F^*})$ (\cite[Theorem 26.3]{Rockafellar}).

$\bullet$ Suppose that $d=1$ and $F$ is convex. Then
\begin{equation} \label{eq: affine dual}
\{v \in \intr(\mathcal{D}_{F^*}): F^* \text{ is affine on } [v-\varepsilon, v+\varepsilon] \text{ for some } \varepsilon >0\} = \intr ( \conv(\im (F'))) \setminus \cl(\im (F')),
\end{equation}
where, recall, $\im(\cdot)$ denotes the image of a function. Thus, kinks of convex functions correspond to affine segments of their convex conjugates, and vice versa. 

In order to prove this, let $v$ belong to the set in the r.h.s.\ of \eqref{eq: affine dual}. Since the function $F'$ is non-decreasing on its domain and $v \notin \cl(\im (F'))$, there is a unique real $u$ such that
\begin{equation}
\label{eq: F' jumps}
\inf (\im (F')) <  F'_-(u) < v < F'_+(u) < \sup (\im (F')).
\end{equation}
Then it is easy to see that $F^*$ is affine on $[F'_-(u), F'_+(u)]$ with slope $u$, which in particular implies that $v$ belongs to the set in the l.h.s.\ of \eqref{eq: affine dual}. For the reverse inclusion, if $v$ belongs to the set in the l.h.s.\ of \eqref{eq: affine dual}, which is open, then $v \in \intr \bigl( \conv(\im (F'))$ by $\mathcal{D}_{F^*} = \cl (\conv(\im (F')))$. 
By taking the Legendre--Fenchel transform of $F^*$ and using \eqref{eq: involution}, which gives $F^{**}=F$ on $\intr(\mathcal D_F)$, we see that \eqref{eq: F' jumps} holds true with $u=(F^*)'(v)$. Hence $v \not \in \cl(\im (F'))$, and thus $v$ belongs to the set in the r.h.s.\ of \eqref{eq: affine dual}.

%Since $K$ is continuously differentiable on $\R^2$, by~Theorem~26.1 and Corollary~26.4.1 in~%\cite{Rockafellar}, we have
%\begin{equation} \label{eq: rint D_I}
%\rint \mathcal{D}_I  = \im(\nabla K).
%\end{equation}
%The set $\rint \mathcal{D}_I$ is convex as the relative interior of a convex set, see~\cite[Theorem~6.2]{Rockafellar}. By \eqref{eq: Taylor K}, the dimension of this set satisfies 
%\begin{equation} \label{eq: dimension}
%\dim(\rint \mathcal{D}_I) = \rank \Sigma = \dim (\supp(X_1)).
%\end{equation} 

%Since $K$ is strictly convex on the affine hull $\aff (\supp(X_1))$,  the implicit function theorem %applied to $\nabla K: \aff (\supp(X_1)) \to \im (\nabla K)$ implies that for any $v \in \im(\nabla %K)$, the equation $\nabla K(u)= v$ has a unique solution $u=u_v \in \aff (\supp(X_1))$ and the %inverse mapping $v \mapsto u_v$ is infinitely differentiable on the relatively open set $\im(\nabla %K)$. Hence
%\begin{equation} \label{eq: I =}
%I(v) = u \cdot v - K(u_v), \qquad v \in \rint \mathcal{D}_I,
%\end{equation}
%and the restriction of $I$ on $\rint \mathcal{D}_I$ is infinitely differentiable.

\subsection{Basic properties of the radial minimum function $\ubar I$} \label{Sec: Basic Proof}

\begin{example}[{\it Discontinuous $\ubar I$}] \label{ex: I_ discontinuous} 
The function $\ubar I$ is not necessarily continuous on $[|\mu|, r_{max})$: it is easy to check that if $\P(X_1=(1,0))=3/4$ and $\P(X_1=(-2,0))=1/4$, then $\ubar I$ has a jump at $r=1$. It is also possible to show that $\ubar I$ is discontinuous for the ``truly'' two-dimensional distribution that is a mixture of the above two-atomic distribution and the uniform distribution on the disk $\{u: |u| \le 1\}$. %We will give a sufficient condition for the continuity of $\ubar I$ in Lemma~\ref{lem: continuity} of the next section. 

It is not clear if $\ubar I$ can be discontinuous for zero mean distributions. 
%Importantly, Example~\ref{ex: I_ non-convex} of Section~\ref{Sec: convexity ubar_I} below gives a zero-mean distribution (supported %on a line) with a {\it non-convex} $\ubar I$.

%For centred distributions, the function $\ubar I$ is continuous (we do not use this statement in the paper and will not prove it) but %not necessarily convex, see Example~\ref{ex: I_ non-convex} in  Section~\ref{Sec: radial proofs} below.
\end{example}

%\begin{example}[{\it Continous but non-convex $\ubar I$}]\label{ex: I_ non-%convex} 
%Put $X_1= (X, 0)$, where $X$ is a random variable with the distribution %belonging to the one-parametric family that satisfies $\P(X \in \{-2,1,3\} ) %=1$, $\E  X = 0$, and $ \E X^3 <0$. The Taylor expansion 
%$$\mathcal{L}((u_1,u_2))=: \mathcal{L}_X(u_1)= 1 + \frac{\E X^2}{2} u_1^2 + %\frac{\E X^3}{6} u_1^3 + o(u_1^3), \qquad u_1 \to 0$$ 
%implies that 
%$\mathcal{L}_X(-u_1) > \mathcal{L}_X(u_1)$ for all $u_1 >0$ that are small %enough enough. Then $I_X(v_1):=I(v_1,0)$ satisfies $I_X(-v_1)< I_X(v_1)$ for %%all $v_1 >0$ small enough; this can be checked directly, and also follows f%rom Proposition 1 by $\bar K(r)=\log [\max(\mathcal{L}_X(-r), \mathcal{L}%_X(r))]$ for $r \ge 0$. Since $\E X^3<0$, we have 
%$$
%I_X(-2) = \log \P(X=-2)> \log \P(X=3)=I_X(3)> I_X(2),
%$$ 
%hence there is a $v_* \in (0,2)$ such that $I_X(v_*)=I_X(-v_*)$. The parameter %can be chosen such that $-I_X'(v_*)\neq I_X(-v_*)$, hence $\ubar I$ is %continuous on $\mathcal D_{\ubar I}$ but non-convex since $(\ubar I)'_-%(v_*)> (\ubar I)'_+(v_*)$.
%\end{example}

\begin{proof}[\bf Proof of Lemma~\ref{lem: properties of I_}.] 
\ref{item: I_ eff domain}) This follows from \eqref{eq: D_I inclusions} and the fact that $\mathcal{D}_{\ubar I}$ is convex.

\ref{item: I_ unimodal}) Clearly, $\ubar I(|\mu|) = 0$ by $I(\mu)=0$ and $I \ge 0$. We claim that for any direction $\ell \in \S^{d-1}$, the function $I_\ell(t):= I(\mu + t\ell)$ is strictly increasing for $t \ge 0$ while it stays finite. 
This implies that $\ubar I$ strictly decreases on $[r_{min}, |\mu|]$ and strictly increases on $[|\mu|, r_{max}]$, since the line segment that joins $\mu$ with a point of minimum of $I$ over the sphere $r \S^{d-1}$  always intersects the sphere $r' \S^{d-1}$ if $0 \le r < r' < |\mu|$ or $|\mu|< r' < r $. 

To prove the claim, we can assume without loss of generality that $\mu =0$ since the rate function of $X_1-\mu$ is $I(v - \mu)$ and the Laplace transforms of $X_1$ and $X_1-\mu$ have the same effective domains.  Since $I_\ell$ is a convex function with minimum at $t=0$, it can cease to be strictly increasing only if it stays zero in a neighbourhood of $0$. If $\ell \cdot X_1 = 0$ a.s., then there is nothing to prove since $I_\ell(t)=+\infty$ for $t>0$ by \eqref{eq: D_I inclusions}, otherwise by the criterion of equality in H\"older's inequality, the function $a \mapsto \log \E e^{a \ell \cdot X_1}$, where $a \in \R$, is strictly convex on its effective domain. Since the interior of this domain contains $0$ by the assumption $0 \in \intr \mathcal D_{\mathcal L}$, it is easy to see that  $I(t\ell) \ge \sup_{a \in \R} (a \ell \cdot t \ell - \log \E e^{a \ell \cdot X_1}) >0$ for $t \neq 0$. This proves the claim.

		\begin{center}
		\includegraphics{fig-poly-6.mps}
		
		\f \label{fig: I_ convex}
	\end{center}	
	
Furthermore, the function $\ubar I$ is convex on $[r_{min}, |\mu|]$ since for any $r_{min} \le r < r' \le |\mu|$ and $\ell \in \ubar{\Lambda}_r, \ell' \in \ubar{\Lambda}_{r'},$ one has
\[
\ubar I(r) + \ubar I(r') = I(r \ell) + I(r' \ell') \ge 2 I \Bigl(\frac{r \ell + r' \ell'}{2} \Bigr) \ge 2 \ubar I \Bigl( \Bigl |\frac{r \ell + r' \ell'}{2} \Bigr | \Bigr) \ge 2 \ubar I  \Bigl( \frac{r+ r' }{2} \Bigr).
\]
Here we used the triangle inequality and the fact that $\ubar I$ decreases on $[0, |\mu|]$, see Figure~\ref{fig: I_ convex} for a geometric explanation in the planar case.  If $\mathcal{D}_{\mathcal L }=\R^d$, then $I$ is strictly convex on its effective domain,  thus the first inequality is strict, hence $\ubar I$ is strictly convex on $[r_{min}, |\mu|]$.

The lower semi-continuity of $\ubar I$ easily follows from that of $I$ using a simple compactness argument. 

\ref{item: I_ discontinuous}) By Part~\ref{item: I_ unimodal}, $\ubar I$ is lower semi-continuous and increasing, and hence left-continuous, on $[|\mu|, r_{max}]$. If $\ubar I$ is discontinuous at an $x \in [|\mu|, r_{max}]$, then it must be $\ubar I(x) < \infty$, otherwise there is a contradiction with the left-continuity of $\ubar I$. For any $\ell \in \ubar{\Lambda}_x$, consider the hyperplane $L$ passing through $x \ell$ and orthogonal to $\ell$. Assume that there is an $v \in L \cap \mathcal D_I$ that is distinct from  $x \ell $. Then $|(1-\varepsilon) x \ell + \varepsilon v|>x$ for every $\varepsilon >0$, hence
\[
\ubar I(x+) \le \lim_{\varepsilon \to 0+} I((1-\varepsilon) x \ell + \varepsilon v) = I(x\ell) = \ubar I(x),
\]
where the first equality holds since the  convex lower semi-continuous function $I$ is continuous on the line segment $[x \ell, v] \subset \mathcal D_I$ (\cite[Corollary~7.5.1]{Rockafellar}). Thus, $\ubar I$ is continuous at $x$, which is a contradiction. Therefore, $L \cap \mathcal D_I = x \ell$, meaning that $x \ell$ is an exposed point of $\mathcal D_I$. It remains to check that $I(x \ell) = -\log \P(X_1 = x \ell)$. 

		\begin{center}
		\includegraphics{fig-poly-702.mps}
		
		\f \label{fig: conv supp}
	\end{center}	

We have $\ell \cdot X_1 \le x$ a.s., where the inequality is strict unless $X_1=x\ell$. Denote by $L_0$ the hyperplane passing through $0$ and parallel to $L$; see Figure~\ref{fig: conv supp}. Let us identify $\R^d$ with $L_0 \oplus \R \ell$. For any $u_1 \in L_0$ such that $\E e^{(u_1 + u_2 \ell ) \cdot X_1} < \infty$ for some real $u_2$, we have
\begin{align} \label{eq: boundary point}
\sup_{u_2 \in \R }  \bigl ( u_2 x - \log \E e^{(u_1 + u_2 \ell ) \cdot X_1} \bigr) 
&=  - \log \bigl ( \inf_{u_2 \in \R }  \E e^{u_1 \cdot X_1 + u_2 (\ell \cdot X_1 - x)} \bigr) \notag\\
&=- \log \E [e^{u_1  \cdot  X_1} \I_{\{X_1 = x \ell\}}] = - \log \P(X_1 = x \ell),
\end{align}
with the second equality following from the dominated convergence theorem using that the random variables in the $\inf \E$ term decrease point-wisely in $u_2 $ since $\ell \cdot X_1 \le x$~a.s. These equalities hold true e.g.\ for $u_1=0$. This yields the required equality
\[
I(x \ell) = \sup_{u_1 \in L_0} \sup_{u_2 \in \R }  \bigl ((u_1 + u_2\ell) \cdot x \ell - \log \E e^{(u_1 + u_2 \ell ) \cdot X_1} \bigr) = - \log \P(X_1 = x \ell),
%\sup_{u_1 \in L_0: \atop (u_1 + \R \ell) \cap \mathcal D_{\mathcal L} \neq %\varnothing} \bigl ( - \log \P(X_1 = x \ell) \bigr),
\]
%hence $I(x \ell) = - \log \P(X = x \ell)$ 
where the last equality holds because no $u_1 \in L_0$ such that $\E e^{(u_1 + u_2 \ell) \cdot X_1} = \infty$ for every $u_2 \in \R$ contributes to the first supremum since for such $u_1$ the l.h.s.\ of the first line in \eqref{eq: boundary point} is $-\infty$. 

Finally, we have $\P(X_1 = x \ell)>0$ by  $x \ell \in \mathcal D_I$.

\ref{item: one element}) Suppose that for an $r \in (r_{min}, |\mu|]$, there are two distinct elements $\ell, \ell'$ in $\ubar{\Lambda}_r$. By convexity of $I$, it holds that $\ubar I(r|\ell + \ell'|/2) \le I(r(\ell + \ell')/2) \le \ubar I(r)$, which is a contradiction by Part~\ref{item: I_ unimodal} since $\ubar I$ is strictly decreasing on $(r_{min}, |\mu|]$ and $|\ell + \ell'|<2$. 
\end{proof} 

\subsection{Radial maxima and minima of conjugate convex functions} \label{Sec: radial proofs}
 
Let us prove the following statement, which may be known in convex analysis  but we found no references. It is stronger than Proposition~\ref{prop: properties of I_} since the Laplace transform of a distribution is lower semi-continuous by Fatou's lemma. In particular, it  applies to distributions with Laplace transform finite only in a neighbourhood of zero.

\begin{prop} \label{prop: generalization}
Let $F: \R^d \to \R \cup \{+ \infty \}$, where $d \ge 1$, be any lower semi-continuous convex function differentiable at $0$ and such that $p_{min}:= \inf\{|u|: F(u)=\infty \}>0$, where $\inf_\varnothing=\infty$ by convention. Put $m:= \nabla F(0)$ and define $\bar F$ and $\ubar{F^*}$ as in~\eqref{eq: radial functions}. Then
\begin{enumerate}[a)]
\item \label{item': conv I_ =} %\label{item: K^- properties} 
$\bar F$ is an increasing convex function on $[0, \infty)$ satisfying $\bar F'_+(0)=|m|$ and 
$$
\conv \ubar {F^*} = (\bar{F})^* \text{ on }[|m|, \infty);
$$
\item \label{item': I_=conv I_} If $r \in \cl (\im(\bar{F}'))$ (and $r \ge |m|$), then $\ubar {F^*}(r) = \conv (\ubar {F^*})(r)<\infty$.
\end{enumerate}

If additionally $F$ is differentiable on $\{u: |u|<p_{min}\}$, and $\bar{\Lambda}_p$ and $\ubar{\Lambda}_r$ are defined for $F$ and $F^*$ as in~\eqref{eq: extreme directions}, then
\begin{enumerate}[a)]
\setcounter{enumi}{2}
\item \label{item': differentiable} For any $p \in (0,p_{min})$, the one-sided derivatives satisfy
$$\bar F_+'(p) = \max_{\ell \in \bar{\Lambda}_p } |\nabla F(p \ell)| \quad \text{and} \quad \bar F_-'(p) = \min_{\ell \in \bar{\Lambda}_p } |\nabla F (p \ell)|.$$
%Importantly, for any $r \ge |\mu|$,
%\item \label{item: conv I_ =} $\conv \ubar I = (\bar{K})^*$ on $[|\mu|, \infty)%$;
%The condition that $\ubar I(r) = \conv \ubar I(r)<\infty$ and $(r, \ubar I(r))$ is an extremal %point of $\epi (\conv \ubar I)$ is equivalent to  $r \in \intr (\im(\bar{K}'))$ (the image is open if we show that $\bar K$ is strictly convex);
\item \label{item': Lambdas} If  $p \in (0,p_{min})$ and $r \ge |m|$ are such that $\bar{F}'(p)=r$, then $\ubar{\Lambda}_r = \bar{\Lambda}_p$.
\end{enumerate}
\end{prop}

The following corollaries to Proposition~\ref{prop: generalization} will easily imply those to   Proposition~\ref{prop: properties of I_}. 

\begin{cor} \label{cor': L differentiable}
If $\bar F$ is differentiable on $(0, p_{min})$, then $\ubar {F^*}$ is strictly convex on $[|m|, (\bar F)'_-(p_{min})]$.
\end{cor}

\begin{cor} \label{cor': condition for diff}
$\bar F$ is differentiable on $(0, p_{min})$ if $F$ is differentiable on $\{u: |u|<p_{min}\}$ and there exists a continuous mapping $\ell: (0, p_{min}) \to \S^{d-1}$ such that $\ell(p) \in \bar \Lambda_p$ for any  $p \in (0,p_{min})$.
\end{cor}

There are few ways to prove Part~\ref{item': conv I_ =} of the proposition, using geometric or analytic approaches. The current simple proof is due to  
Fedor Petrov.
 
\begin{proof}[\bf Proof of Proposition~\ref{prop: generalization}.]
\ref{item': conv I_ =}) We have $\bar F(p) = \sup_{\ell \in \S^{d-1}} F(p \ell)$, where $p \ge 0,$ hence $\bar F$ is convex as a maximum of convex functions $F_\ell(\cdot)  := F(\cdot \ell)$. Furthermore, the convex function $F$ attains its maximum over any closed compact convex set on the boundary of the set. Therefore for any $0 \le p < p' \le p_{min}$, we have 
	\[
	\bar F(p)=\max_{u: |u|=p }F(u)=\max_{u: |u|\le p}F(u) \le \sup_{u: |u|\le p'}F(u)=\bar F(p')\]
(where the supremum may not be attained if $p'=p_{min}$). 
	Hence $\bar F$ is increasing on $[0, \infty)$ since $\bar F(p)=\infty$ for $p> p_{min}$. The right derivative of $\bar F$ at $0$ clearly satisfies $\bar F'_+(0) = |m|$.
	
It remains to prove that $\conv(\ubar{F^*})= (\bar F)^*$ on $[|m|, \infty)$. We first claim that
\begin{equation} 
\label{eq: Fedia}
(\ubar{F^*})^*(p) = \bar F(p), \qquad p \ge 0.
\end{equation} 
In fact, by the definition, we have $\ubar{F^*} (r)= \infty$  for $r <0$, hence for every real $p$,
\begin{align*}
(\ubar{F^*})^*(p) &= \sup_{r \ge 0} \bigl (pr - \ubar{F^*}(r) \bigr) = \sup_{r \ge 0} \bigl (pr - \inf_{\ell \in \S^{d-1}} F^*(r \ell) \bigr) \\
&= \sup_{r \ge 0, \ell \in \S^{d-1}} \bigl (pr - F^*(r \ell) \bigr) = \sup_{v \in \R^d} \bigl (p |v| - F^*(v) \bigr). 
\end{align*}
On the other hand, by the assumptions, $F$ is convex, lower semi-continuous, and has non-empty effective domain,  hence  $F=F^{**}$ holds  by \eqref{eq: involution}. Then \eqref{eq: Fedia} follows since for $p \ge 0$, 
\[
\bar F(p) = \sup_{\ell \in \S^{d-1}} F(p \ell) = \sup_{\ell \in \S^{d-1}} \sup_{v \in \R^d} \bigl( p \ell \cdot v -F^*(v)  \bigr)= \sup_{v \in \R^d} \bigl (p |v| - F^*(v) \bigr).
\]

By the definition, we have $\bar F(p)=\infty$  for $p<0$, 
hence the Legendre--Fenchel transform of $\bar F$ is fully defined by the values of $\bar F$ on $[0, \infty)$. Likewise, the Legendre--Fenchel transform of $(\ubar{F^*})^*$ restricted to $[|m|, \infty)$ is  defined by the values of $(\ubar{F^*})^*$ on $[0, \infty)$. In fact, the function $p \mapsto pr - (\ubar{F^*})^*(p)$ is increasing on $(-\infty, 0]$ for any $r \ge |m|$ since $(\ubar{F^*})^*$ is a convex function, whose right derivative increases on its domain and its value at $0$ equals that of $\bar F$ by \eqref{eq: Fedia}, while we already proved that $\bar F'_+(0)=|m|$. Therefore, \eqref{eq: Fedia} implies $(\ubar{F^*})^{**} = (\bar F)^*$ on $[|m|, \infty)$, hence $\cl (\conv(\ubar{F^*}))= (\bar F)^*$ on $[|m|, \infty)$ by \eqref{eq: involution}.

It remains to remove the closure operation $\cl$ from the last equality. It suffices to show that $\conv(\ubar{F^*})$ is lower semi-continuous at the boundary points of its effective domain, which coincides with that of $\ubar{F^*}$. By a simple compactness argument it follows from the lower semi-continuity of $F^*$ that $\ubar{F^*}$ is also lower semi-continuous. Every point $x \in \partial (\mathcal{D}_{\ubar {F^*}})$ has a neighbourhood $U$ such that the convex function $\conv(\ubar{F^*})$ is either strictly increasing, strictly decreasing, or constant on $U \cap \mathcal{D}_{\ubar {F^*}}$. From the definition of the largest convex minorant, it follows that $\conv(\ubar{F^*})(x)= \ubar{F^*}(x)$. This equality, combined with the property of lower semi-continuity of $\ubar{F^*}$ at $x$, implies the same property for $\conv(\ubar{F^*})$ by a simple consideration of the three cases mentioned above.
	
\ref{item': I_=conv I_})  It it easy to see that the non-negative function $\ubar {F^*} - \conv (\ubar {F^*})$, which we define to be zero outside $\mathcal{D}_{\ubar {F^*}}$, is lower semi-continuous. In fact, this property holds at the points of $\intr (\mathcal{D}_{\ubar {F^*}})$ by continuity of $\conv (\ubar {F^*})$ and lower semi-continuity of $\ubar {F^*}$, which we showed above in the proof of Part~\ref{item': conv I_ =}. At the points of $\partial (\mathcal{D}_{\ubar {F^*}})$, this is true by non-negativity and the fact that $\ubar {F^*} = \conv (\ubar {F^*})$ on $\partial (\mathcal{D}_{\ubar {F^*}})$. Hence the set 
\[
\{r \ge |m|: \ubar {F^*}(r) = \conv (\ubar {F^*})(r)\} = \{r \ge |m|: \ubar {F^*}(r) -\conv (\ubar {F^*})(r) \le 0\}
\]
is closed as a sub-level set of a lower semi-continuous function. Therefore, if $\ubar {F^*}(r) >\conv (\ubar {F^*})(r)$ for an $r  \in (|m|, \sup \mathcal{D}_{\ubar {F^*}})$, then this inequality also holds on an open interval $(r_1, r_2)$  that contains $r$, on which $\conv (\ubar {F^*})$ must be affine. Since $\conv (\ubar {F^*})=(\bar F)^*$ on $[|m|, \infty)$ by Part~\ref{item': conv I_ =}, we conclude that $(\bar F)^*$ is affine on $[r_1, r_2]$. As we explained in Section~\ref{Sec Convex basics}, this yields that $r \not \in \cl(\im( \bar F'))$, which is a contradiction.

\ref{item': differentiable}) The one-sided derivatives of $\bar F$ exist by convexity of this function proven in Part~\ref{item': conv I_ =}. The set $\bar{\Lambda}_p:= \argmax_{\ell \in \S^{d-1}} F(p \ell)$ is well-defined
since the convex function $F$ is continuous on $\rint( \mathcal{D}_F) $ and $p \S^{d-1} $ is a compact subset of $\rint( \mathcal{D}_F) $ by  $p \in (0, p_{min})$.

For any $\ell \in \bar{\Lambda}_p$, the gradient $\nabla F(p \ell)$ is directed along $\ell$ since $p \ell$ is an extremal point of the function $F$ over the sphere $p \S^{d-1}$ and $F$ is differentiable, and hence continuously differentiable, on $\{u: |u|<p_{min}\}$; see Section~\ref{Sec Convex basics}. Hence $|\nabla F(p \ell)| = F_\ell'(p)$ and by
\[
\bar F'_+(p) = \lim_{\varepsilon \to 0+} \varepsilon^{-1}(\bar F(p+\varepsilon)- \bar F(p)) \ge \lim_{\varepsilon \to 0+} \varepsilon^{-1}(F((p+\varepsilon)\ell)- F(p\ell)) =  F_\ell'(p),
\]
we arrive at $\bar F'_+(p) \ge \max_{\ell \in \bar{\Lambda}_p} |\nabla F(p \ell)|$, where the r.h.s.\ accounts the fact that the function $\nabla F$, which is continuous on $p \S^{d-1}$, attains its maximum on the compact set $p \bar{\Lambda}_p$.

Furthermore, since $\S^{d-1}$ is compact, there exist two sequences $p_k \to p+$ and $\ell(k) \in \bar{\Lambda}_{p_k}$ such that $\ell(k) \to \ell$ for some $\ell \in \S^{d-1}$ as $k \to \infty$. Then necessarily $\ell \in \bar{\Lambda}_p$ since $F$ and $\bar F$ are continuous on some neighbourhoods of $p \S^{d-1}$ and $p$, respectively. Finally,
\[
\bar F(p_k) - \bar F(p)= F(p_k \ell(k)) - F(p \ell) = (p_k \ell(k) - p \ell) \cdot (\nabla F(p \ell) + o(1)) \le (p_k - p) (|\nabla F(p\ell)| + o(1))
\]
as $k \to \infty$, and thus $\bar F'_+(p) \le \max_{\ell \in \bar{\Lambda}_p} |\nabla F(p \ell)|$. This inequality, combined with the opposite one proven above, yields the equality required.

The argument for $\bar F'_-(p)$ is analogous. 

\ref{item': Lambdas}) The set  $\ubar{\Lambda}_r :=\argmin_{\ell \in \S^{d-1}} F^*(r \ell)$ is well-defined since $F^*$ is lower semi-continuous. First check  that $\bar \Lambda_p \subset \ubar {\Lambda}_r$. For any $\ell \in \bar \Lambda_p$, $\nabla F(p\ell)$ is directed along $\ell$, hence by Part~\ref{item': differentiable} it holds that $\nabla F(p\ell)=r\ell$. 
Note that $F(u)\geq F(p\ell)+r\ell\cdot(u-p\ell)$ for any $u \in \mathbb{R}^d$ since the right-hand side of this inequality defines the support hyperplane to graph of $F$ at the point $(p\ell, F(p\ell))$. Then
\begin{align}
	\ubar{ F^*}(r) \le F^*(r\ell)=\sup_{u\in \mathbb{R}^d} \bigl(  r\ell\cdot u - F(u) \bigr) &\leq
	\sup_{u\in \mathbb{R}^d} \bigl(  r\ell\cdot u -F(p\ell)-r\ell\cdot(u-p\ell) \bigr) \notag \\
	&= r\ell\cdot p\ell-F(p\ell)=rp-\bar F(p). \label{eq: technical}
\end{align}
The concave function $q \mapsto r q -\bar F(q) $ attains its maximum at $q=p$ since by the assumption, it holds that $\bar F'(p)=r$. Then by Parts~\ref{item': conv I_ =} and~\ref{item': I_=conv I_},
\begin{equation} \label{eq: aux}
rp-\bar F(p) = (\bar F)^*(r)=\conv (\ubar{F^*})(r) = \ubar{F^*}(r),
\end{equation}
and since the latter expression equals the first term in~\eqref{eq: technical}, we get $\ell \in \ubar {\Lambda}_r(F^*)$.

It remains to prove the reverse inclusion $\ubar {\Lambda}_r \subset \bar \Lambda_p $. Suppose that $\ell \in \ubar {\Lambda}_r$. Combining the Fenchel inequality with \eqref{eq: aux}, we obtain
\begin{equation*}
	\bar F(p) \ge  F(p\ell) \geq
	r\ell\cdot p\ell -F^*(r\ell)
	= rp - \ubar {F^*}(r) = \bar F(p),
\end{equation*}
which implies that $\ell \in \bar \Lambda_p$.
\end{proof}

\begin{proof}[{\bf Proof of Corollary~\ref{cor': L differentiable}.}]
Since the function $\bar F$, which is convex on $[0, \infty)$ by Proposition~\ref{prop: generalization}.\ref{item': conv I_ =}, is assumed to be differentiable on $(0, p_{min})$,  it is continuously differentiable there; see Section~\ref{Sec Convex basics}. Then by \eqref{eq: affine dual}, $(\bar F)^*$  is affine on no subinterval of 
\[
[\inf(\im(\bar F')), \sup(\im(\bar F'))] = [(\bar F')_+(0), (\bar F')_-(p_{min})]=[|m|, (\bar F')_-(p_{min})],
\]
and therefore strictly convex there. So is the function $\ubar {F^*}$, which equals $(\bar F)^*$ on $[|m|, \infty)$ by Proposition~\ref{prop: generalization}.\ref{item': conv I_ =} and~\ref{prop: generalization}.\ref{item': I_=conv I_}. 
\end{proof}

\begin{proof}[{\bf Proof of Corollary~\ref{cor': condition for diff}.}]
Since the function $\bar F$ is convex on $[0, \infty)$, its left and right derivatives satisfy (\cite[Theorem~24.1]{Rockafellar})
\[
\bar{F}'_+(p-) = \bar{F}_-'(p) \le \bar{F}_+'(p) = \bar{F}_-'(p+), \qquad p \in (0, p_{min}).
\]
On the other hand, we have $\bar{F}_-'(p) \le   |\nabla F(p \ell(p))| \le \bar{F}_+'(p)$ by Proposition~\ref{prop: generalization}.\ref{item': differentiable}. The claim follows by combining these relations and using that $|\nabla F(p \ell(p))|$ is continuous on $(0, p_{min})$, which is true since $p\ell(p)$ is continuous on $(0, p_{min})$ and $\nabla F$ is continuous on $\intr(\mathcal D_F)$; see Section~\ref{Sec Convex basics}.
\end{proof}

\begin{proof}[{\bf Proofs of Corollaries~\ref{cor: L differentiable} and~\ref{cor: condition for diff}.}]
We apply Corollaries~\ref{cor': L differentiable} an~\ref{cor': condition for diff} with $K$ substituted for $F$. Since $\mathcal D_{\mathcal L}=\R^d$ by  the assumption, we have $p_{min}=\infty$ by $K=\log \mathcal L$. Then Corollary~\ref{cor: condition for diff} follows from  Corollary~\ref{cor': condition for diff}. Furthermore, it follows from Proposition~\ref{prop: properties of I_}.\ref{item: differentiable} that $\lim_{p \to \infty}(\bar K)'_-(p)=r_{max}$. Then $\ubar I$ is strictly convex on $[|\mu|,r_{max}]$ by Corollary~\ref{cor': L differentiable}, while $\ubar I$ is strictly convex on $[r_{min}, |\mu|]$ by
Lemma~\ref{lem: properties of I_}.\ref{item: I_ unimodal}. Since $I$ attains its minimum at $\mu$ and is continuous at $\mu$ except for the trivial case $X_1 = \mu$ a.s.,  $\ubar I$ is strictly convex on the interval $[r_{min}, r_{max}]$, which contains $\mathcal D_{\ubar I}$ by Lemma~\ref{lem: properties of I_}.\ref{item: I_ eff domain}. This proves Corollary~\ref{cor: L differentiable}, which is trivial in the remaining case $X_1 = \mu$~a.s.
\end{proof}

\subsection{Convexity of the radial minimum function $\ubar I$} \label{Sec: convexity ubar_I}

Here we prove that $\ubar I$ is convex for the distributions described in Proposition~\ref{prop: conjecture holds}.

\begin{proof}[{\bf Proof of Proposition~\ref{prop: conjecture holds}.}]

\ref{item: pseudo-affine}) Denote $r:=\rank A$. We assume that $r \ge 1$, otherwise the claim is trivial. The $d \times k$ matrix $A$ admits a singular value decomposition $A = U D V^\top$, where $D$ is an $r \times r$ diagonal matrix whose diagonal entries are non-zero singular values of~$A$ (i.e., the square roots of non-zero eigenvalues of $A A^\top$), and $U$ is $d \times r$ matrix and $V$ is a $k \times r$ matrix such that both $U^\top U$ and $V^\top V$ are the unit $r \times r$ matrices. 

Put $L:= U \R^r$. Then the assumption $A A^\top \mu =\|A A^\top \! \| \mu$ implies that $\mu \in L$. Furthermore, it is easy to check that $(U^\top u_1, U^\top u_2) = (u_1, u_2)$ for any $u_1, u_2 \in L$ and $|U^\top u|<|u|$ for $u \in \R^d \setminus L$.  Therefore, from the equalities
\[
\|A A^\top \! \| = \max_{u \in \S^{d-1}} |A A^\top u| =  \max_{u \in \S^{d-1}} |U D^2 U^\top u| =   \max_{u \in \S^{d-1} \cap L} |U D^2 U^\top u| =
\max_{u' \in \S^{r-1}} |D^2 u'| = \| D\|^2,
\]
and  the assumption $A A^\top \mu =\|A A^\top \! \| \mu$, we see that the vector $\mu':= U^\top \mu$ in $\R^r$ satisfies $D \mu' = \sigma_1 \mu'$, where $\sigma_1:=\|D\|$ is the largest singular value of $A$.

Then for any $v \in L$, by $X_1 \in L$ a.s.\ we have
\begin{align*}
I(v)= \sup_{u \in \R^d} \Big( u \cdot v - \log \E e^{u \cdot X_1} \Big) 
&=\sup_{u \in L} \Big( u \cdot (v-\mu) - \log \E e^{u \cdot A Y_1} \Big) \\
&= \sup_{ u \in L} \Big( U^\top u \cdot  U^\top (v - \mu) - \log \E e^{DU^\top u \cdot V^\top Y_1} \Big). 
\end{align*}
Denote by $J$ the rate function of the random vector $V^\top Y_1$ in $\R^r$.
Let us use that $I(v)=+\infty$ for $v \in \R^d \setminus L$ (by $X_1 \in L$ a.s.) and change variables $u' = U^\top u$, $v'= U^\top v$, $u''=Du'$ to get
\begin{align*}
\ubar I (r) &= \min_{v \in r \S^{d-1} \cap L} I(v) = 
\min_{v' \in r \S^{r-1}}  \sup_{u' \in \R^r} \Big( u' \cdot (v' -  \mu') - \log \E e^{D u' \cdot V^\top Y_1} \Big)  \\
&  = \min_{v' \in r \S^{r-1}}  \sup_{u'' \in \R^r} \Big( D^{-1} u'' \cdot (v' -  \mu') - \log \E e^{u''\cdot V^\top Y_1} \Big) = \min_{v' \in r \S^{r-1}} J \big(D^{-1} (v' - \mu')\big ),
\end{align*}
where in the last equality we also used that $D^{-1}$ is symmetric. 

The distribution of  $V^\top Y_1$ on $\R^r$ is rotationally invariant since so is that of $Y_1$ on $\R^k$. Therefore, $J$ is a radial function, hence
\[
\ubar I (r) =  \ubar J \Big( \min_{v' \in r \S^{r-1}}  |D^{-1} (v' -  \mu')|\Big ) = \ubar J \big(\sigma_1^{-1}|r - |\mu'||\big), \qquad r \ge 0,
\]
where  $|\mu'|=|\mu|$ and we used that $D^{-1} \mu' = \sigma_1^{-1} \mu'$ and $\sigma_1^{-1}$ is the smallest eigenvalue of $D^{-1}$. Hence $\ubar I$ is convex since so is $\ubar J$ and $\ubar J$ is increasing on $[0, \infty)$.

\ref{item: Gaussian})  For $d=1$, $\ubar I$ is convex by Part~\ref{item: pseudo-affine}, so we assume that $d \ge 2$.  We will give a detailed treatment for illustrative purposes in the planar case  and then proceed to higher dimensions.

1. The planar case $d=2$ with non-degenerate covariance matrix $\Sigma$ of $X_1$. 

The cumulant generating function $K$ of a Gaussian$(\mu, \Sigma)$ distribution is $K(u) = u^\top \mu + \frac12 u^\top \Sigma u$. 
%The graph of this function is an elliptic paraboloid. 
By Corollary~\ref{cor: L differentiable} , it suffices to show that the radial maximum function $\bar K(p)$ is differentiable on $(0, \infty)$. Since $\bar K$ is invariant under orthogonal transformations of $\R^2$, without loss of generality we can assume that 
\[
K(x,y)=\frac12 a(x-x_0)^2+ \frac12 b(y-y_0)^2 + c,
\]
where $a$ and $b$ are the eigenvalues of $\Sigma$, $c=-\frac12 a x_0^2 - \frac12 b y_0^2$, $\mu = (-ax_0, -by_0)$, and $x_0, y_0 \ge 0$ by changing directions of the axes, if necessary. We can further assume that $a > b>0$ and $x_0 + y_0>0$, since the cases $a=b$ and $x_0=y_0=0$ are already covered by Part~\ref{item: pseudo-affine}.

To prove that $\bar K$ is differentiable, by Corollary~\ref{cor: condition for diff} it suffices to show that there is a continuous path $\ell(p)$ on the unit sphere that belongs to $\bar \Lambda_p $ for every $p>0$. Suppose that $\ell \in \bar \Lambda_p$, i.e.\ $K$ attains its maximum over $p \S^1$ at the point $p \ell$. Then $\nabla K(p \ell) = t p\ell$ for some non-zero real $t$, that is $(a(x-x_0),b(y-y_0)) = (t x,ty)$. Equivalently, 
	\begin{equation}
		\label{eq:hyperbola polynom}
            (a-t) x = a x_0, \quad (b-t) y = b y_0.
	\end{equation}

a)  The case $x_0, y_0 >0$. The set $\bar \Lambda_p$ lies in the quadrant $\{(x, y):x\leq 0, y\leq 0\}$ since $K(-|x|,-|y|) < K(x,y)$ for any pair $(x, y)$ in the complement of the quadrant. Hence, because the right-hand sides of the equalities in \eqref{eq:hyperbola polynom} are strictly positive, we have $t>a$. Therefore, equalities \eqref{eq:hyperbola polynom} define the curve 
	\begin{equation}
		\label{eq:hyperbola parametric}
            h(t):=  \Big( \frac{a x_0}{a-t}, \frac{b y_0}{b-t} \Big), \qquad t>a,
	\end{equation}
marked in bold in Figure~\ref{fig: Apollonian hyperbola}. Note in passing that $h(t)$ is an arc of the Apollonian hyperbola for the ellipses that are contour lines of $K$; see Glaeser et al.~\cite[Section 9.3]{glaeser2016universe} for details.

\begin{center}
	\parbox{7cm}{
	\begin{center}
	\includegraphics{fig-poly-12.mps}\\
	\f \label{fig: Apollonian hyperbola}
	\end{center}
	}
	\parbox{7cm}{
	\begin{center}
	\includegraphics{fig-poly-13.mps}\\

	\f \label{fig: degenerate Apollonian}
	\end{center}	
	}
\end{center}

Both coordinates of $h(t)$ are strictly decreasing and continuous in $t$, hence the function $t \mapsto |h(t)|$ is a strictly decreasing continuous bijection from $(a, \infty)$ to $(0, \infty)$. Therefore, the curve in \eqref{eq:hyperbola parametric}, obtained from the necessary condition~\eqref{eq:hyperbola polynom} for a maximum, has a unique point of intersection with the circle $p \S^1$. This point must be the unique element of the non-empty set $\bar \Lambda_p $. Thus, $\ell(p):=h(|h|^{-1}(p))$ is the curve required.

b) The cases $x_0>0, y_0=0$ and $x_0=0, y_0>0$.  We  consider them solely for the purpose of illustration since they will be covered below in Part 2 using a different general argument. Meanwhile, note in passing that here equalities \eqref{eq:hyperbola polynom} define two lines $x=\frac{ax_0}{a-b}$ and $y=-\frac{by_0}{a-b}$, which can be regarded as the limit shapes for the hyperbolas in \eqref{eq:hyperbola parametric}.

%Here the hyperbola $\mathcal{H}$ degenerates into two lines given in~%\eqref{eq:axis of hyperbola}. 
%Below we will present a different method which covers such generate cases %in any dimension so the following argument is given only for completeness %of consideration. 
	
It is easy to see that in the first case $y_0=0$, we have $\bar \Lambda_p = \{(-1, 0)\}$ for every $p>0$, so $\ell(p):=(-1,0)$; this situation is actually covered above in Part~\ref{item: pseudo-affine} . In the second case $x_0=0$, from \eqref{eq:hyperbola polynom} we have $x=0$ or $a=t$. Both solutions contribute to the answer -- we have $\bar \Lambda_p = \{(0,-1)\}$ for $p \in (0, b y_0/(a-b)]$, and $\bar \Lambda_p$ consists of two directions symmetric about the $y$-axis for $p >b y_0/(a-b)$. The set $\cup_{p>0} p\bar \Lambda_p$ is marked in bold in Figure~\ref{fig: degenerate Apollonian}. 
Clearly, there is a continuous path $\ell(p)$ of directions such that  $\ell(p) \in \bar \Lambda_p$ for every $p>0$, as required. 

2. Arbitrary dimension $d \ge 2$ with non-degenerate $\Sigma$.

Take a basis of $\R^d$ of eigenvectors of $\Sigma$, where  the coordinates $\mu_i$ of $\mu$ are non-positive. If all $\mu_i$'s are strictly negative, we argue exactly as above, putting $\ell(p):=h(|h|^{-1}(p))$, where
\[
h(t):=  -\Big( \frac{\mu_1}{\sigma_1^2-t}, \ldots, \frac{\mu_d}{\sigma_d^2-t} \Big), \qquad t>\sigma_1^2,
\]
and $\sigma_1^2 \ge \ldots \ge \sigma_d^2 >0$ are the eigenvalues of $\Sigma$. Note in passing that if $\sigma_1^2 = \ldots =\sigma_d^2 $, then $h(t)$ parametrizes the half-line emanating from $0$ in the direction of $\mu$.

If some coordinates of $\mu$ are zero, we proceed differently from our consideration in the planar case and prove the convexity of $\ubar I$ directly. 
%rather than using Corollary~\ref{cor: condition for diff}. 
The rate function of the Gaussian$(\mu, \Sigma)$ distribution is given by $I(v)=\frac12 (v-\mu)^\top \Sigma^{-1}(v-\mu)$. For any $\varepsilon >0$, the Gaussian$(\mu-\varepsilon e_d, \Sigma)$ distribution, where $e_d:=(1, \ldots, 1)$, has the rate function $I_\varepsilon(v):=I(v+\varepsilon e_d)$. All coordinates of $\mu-\varepsilon e_d$ are strictly negative, hence each function $\ubar{I_\varepsilon}$ is convex on $[0, \infty)$ as shown above. On the other hand, $I_\varepsilon \to I$ as $\varepsilon \to 0+$ uniformly on every compact subset of $\R^d$ since $I$ is continuous on $\R^d$, and hence locally uniformly continuous. Then $\ubar{I_\varepsilon}(r) \to \ubar I(r)$ for every $r  \ge 0$, which implies that $\ubar I$ is  convex on $[0, \infty)$, as required. 

3. Arbitrary dimension $d \ge 2$ with degenerate $\Sigma$.

Put $L:=\Sigma \R^d$ and note that $\Sigma_L:={\Sigma|}_L$ is a bijection from $L$ to $L$. The rate function of the Gaussian$(\mu, \Sigma)$ distribution with degenerate $\Sigma$ is given by $I(v)=\frac12 (v-\mu)^\top \Sigma_L^{-1}(v-\mu)$ for $v \in \mu + L$ and $I(v)=+\infty$ for $v \not \in \mu + L$.

For any $\varepsilon >0$, let $\Sigma_\varepsilon$ be the positive definite $d \times d$ matrix defined by $\Sigma_\varepsilon u = \Sigma u$ for $u \in L$ and $\Sigma_\varepsilon u=\varepsilon u$ for $u \in \ker \Sigma$. Let $I_\varepsilon$ be the rate function of the Gaussian$(\mu, \Sigma_\varepsilon)$ distribution. We have $I_\varepsilon(v)=I(v)$ for $v \in \mu +L$ and $I_\varepsilon(v) \nearrow I(v)$  as $\varepsilon \to 0+$ for $v \not \in \mu +L$. 
Since the matrix $\Sigma_\varepsilon$ is non-degenerate, each function $\ubar{I_\varepsilon}$ is convex on $[0, \infty)$ as shown above.  To conclude that $\ubar I$ is convex on $[0, \infty)$, it remains to prove that $\ubar{I_\varepsilon}(r) \to \ubar I(r)$ for every $r  \ge 0$. 

Denote by $v'$ the orthogonal projection of a $v \in \R^d$ on $L$ and put $v'':=v-v'$. Then 
\[
I_\varepsilon(v)=\frac12 (v'-\mu')^\top \Sigma_\varepsilon^{-1}(v'-\mu') + \frac12 \varepsilon^{-1} |v'' - \mu''|^2 = I(v'+\mu'') + \frac12 \varepsilon^{-1} |v'' - \mu''|^2.
\]
Fix an $r > 0$. 
Then for all $\varepsilon >0$ small enough, we have 
\[
\ubar{I_\varepsilon}(r) = \min_{v \in r \S^{d-1} } I_\varepsilon(v) \ge \min_{\substack{v \in r \S^{d-1}:  \\ |v'' - \mu''|\le \varepsilon^{1/3}}  }  I(v'+\mu'') \ge \min_{\substack{v \in r \S^{d-1}:  \\ |v'' - \mu''|\le \varepsilon^{1/3}}  }  \ubar I \Big(\sqrt{|v'|^2+|\mu''|^2} \Big),
\]
where in the first equality we used no $v$ such that  $|v'' - \mu''| > \varepsilon^{1/3}$ contributes to the first minimum since 
$I_\varepsilon(v)>\frac12  \varepsilon^{-1/3}$ for such $v$. Finally, since $|v'|^2=r^2 - |v''|^2$ for $v \in r \S^{d-1}$, 
\[
\ubar{I_\varepsilon}(r) \ge \min_{\substack{v \in r \S^{d-1}:  \\ |v'' - \mu''|\le \varepsilon^{1/3}}  }  \ubar I \Big(\sqrt{r^2+|\mu''|^2-|v''|^2} \Big) \ge \min_{|\delta|\le 2 \varepsilon^{1/3}  |\mu''| + \varepsilon^{2/3} }  \ubar I \big(\sqrt{r^2+\delta} \big). 
\]

Then, since $\ubar I(r) \ge \ubar{I_\varepsilon}(r) $ and $\ubar I$ is continuous at $r$, we obtain that $\ubar{I_\varepsilon}(r) \to \ubar I(r)$ as $\varepsilon \to 0+$ for every $r>0$, as required. This is also true for $r=0$ since $I_\varepsilon(0) \to I(0)$.
\end{proof}

\section{Proofs of the main results} \label{Sec: proof}
\subsection{Basic facts on large deviations} \label{Sec: basics on LDP} ~

$\bullet$ Let $(Z_n)_{n \ge 1}$ be random elements of a Polish space $\mathcal M$ equipped with a metric $d$, and let $\mathcal{J}:\mathcal M\to[0,\infty]$ be a lower semi-continuous function. We say that $\mathcal J$ is  {\it tight} if  its sub-level sets $\{x\in \mathcal M \,:\, \mathcal{J}(x) \leq \alpha \}_{\alpha\ge 0}$ are compact. We say that the collection $(Z_n)_{n \ge 1}$ satisfies a {\it large deviations principle (LDP,}  in short) in $\mathcal M$ with speed $n$ and the {\it rate function} $\mathcal{J}$ if for every Borel set $B\subset \mathcal M$,
\begin{align}\label{eq: LD general}
-\inf_{x\in \intr B}\mathcal{J}(x) \leq\liminf_{n\to\infty}{1\over n}\log \P(Z_n\in B) & \leq\limsup_{n\to\infty}{1\over n}\log \P(Z_n\in B) \leq-\inf_{x\in \cl B}\mathcal{J}(x),
\end{align}
where, as usual, we agree that $\inf_\varnothing = +\infty$. 
%Likewise, $(Z_n)_{n \ge 1}$ satisfies a {\it weak LDP} if the lower bound in \eqref{eq: LD general} holds for every $B$ and the upper bound is satisfied only if the set $\cl B$ is compact. 
We assume throughout that $\mathcal J$ is tight; so are all the rate functions considered in this paper. A Borel set $B \subset \mathcal M$ is called {\it regular} for the rate function $\mathcal J$ if the infima in \eqref{eq: LD general} are equal. Since $\mathcal{J}$ is tight, the infimum on the r.h.s.\ of \eqref{eq: LD general} is always attained at some~$x$.

$\bullet$ Assume that $B \subset \mathcal M$ is a closed set such that $\lim_{n\to\infty} \frac1n \log \P\big(Z_n\in B\big) =-\inf_{x\in B}\mathcal{J}(x)$ (e.g., we can take any regular closed set) and $B \cap \mathcal D_{\mathcal{J}} \neq \varnothing$. Then for any $\varepsilon >0$,
\begin{equation} \label{eq: LDP conditional}
\lim_{n\to\infty}{1\over n}\log \P\Big(d(Z_n,x) \le \varepsilon \text{ for some } x \in B \text{ such that } \mathcal{J}(x) = \min_{y \in B} \mathcal{J}(y) \Big | \Big. Z_n \in B \Big) =1. 
\end{equation}
This means that given the large deviations event $\{Z_n \in B\}$, the random elements $Z_n$ asymptotically concentrate around the compact set of minimizers of the rate function $\mathcal{J}$ over $B$. This follows from \eqref{eq: LD general} since the conditioned event in \eqref{eq: LDP conditional} is $\{d(Z_n, \argmin_{x \in B} \mathcal{J}(x)) \le\varepsilon \big\}$ and we have
\[
\min_{x \in B} \mathcal{J}(x) < \inf_{x \in B : \, d(x, \, \argmin_{y \in B} \mathcal{J}(y) ) \ge \varepsilon} \mathcal{J}(x).
\]
The last inequality holds true since by tightness of $ \mathcal{J}$, the infimum  on the r.h.s.\ is attained on some $x \not \in \argmin_{y \in B} \mathcal{J}(y)$ unless the minimum is taken over the empty set, in which case the r.h.s.\ is $+\infty$ and the inequality is still true.

$\bullet$  Denote by $C_0[0,1]=C_0([0,1];\R^2)$ the space of continuous functions $h:[0,1] \to \R^2$, i.e.\ planar curves, that satisfy $h(0)=0$. We equip this space with the usual metric of uniform convergence. Denote by  $AC_0[0,1]$ its subspace  of functions with absolutely continuous coordinates. Let $S_n(\cdot) \in C_0[0,1]$ be the random piecewise linear functions that satisfy $S_n(k/n):=S_k$, $0 \le k \le n$, where $S_0:=0$, and their values at the other points of $[0,1]$ are defined by linear interpolation. Define the function $I_C: C_0[0,1] \to [0, \infty]$ to be
\begin{equation} \label{eq: I def}
I_C(h):=\left\{
                 \begin{array}{ll}
                  \int_0^1 I(h'(t)) dt, & \text{if } h \in AC_0[0,1]; \\
                  +\infty, & \text{otherwise}.
                 \end{array}
               \right.
\end{equation}

The following result, although stated in a different form, is due to Mogulskii~\cite[Theorem 2, Part II]{Mogulskii}; there were earlier works in this direction by A.A.\ Borovkov. The exact statement presented here appears in the proof of Theorem 5.1.2 in book by Dembo and Zeitouni~\cite{DemboZeitouni}. 

\begin{theorem*}[Mogulskii's LDP]
Assume that $X_1$ is a random vector in $\R^d$, $d \ge 1$, such that $\mathcal D_{\mathcal L} = \R^d$. Then the sequence of random functions $(S_n(\cdot)/n)_{n \ge 1}$ satisfies the LDP in $C_0[0,1]$ with speed $n$ and the tight rate function $I_C$.
\end{theorem*}

$\bullet$ The above LDP for the trajectories $S_n(\cdot)/n$ readily implies that the random vectors $(S_n/n)_{n \ge 1}$ satisfy the LDP in $\R^2$ with speed $n$ and the tight rate function $I_1(v) := \inf_{h: h(1)=v} I_C(h)$ for $v \in \R^2$. This follows by applying  the {\it contraction principle}  (\cite[Theorem~4.2.1]{DemboZeitouni}) and continuity of the mapping $h \mapsto h(1)$. Then $I_1=I$ by Jensen's inequality: 
\begin{equation} \label{eqn: Jensen direct}
I_C(h)=\int_0^1 I(h'(t))dt \ge I(h(1)), \qquad h \in AC_0[0,1]. 
\end{equation}
In particular, for any Borel set $B \subset \R^2$ that is regular for the rate function $I$, we have
\begin{equation*} %\label{eq: LD regular R^2}
\lim_{n \to \infty} \frac{1}{n} \log \P(S_n/n  \in B)  =- \inf_{v \in B} I(v). 
\end{equation*} 
Since the rate function $I$ is strictly convex, by \eqref{eq: LDP conditional} this implies that the trajectories $S_n(\cdot)$ that result in the large deviations event $\{S_n/n \in B\}$ are asympotically linear, as in \eqref{eq: LD shape <}.

%\begin{remark} \label{rem: convex regular}
%By Cram{\'e}r's theorem (\cite[Theorem~6.13 and Corollary~6.1.6]{DemboZeitouni}), \eqref{eq: LD regular R^2} actually %holds merely under $0 \in \intr(\mathcal D_{\mathcal L})$ rather than $\mathcal D_{\mathcal L} = \R^2$. Moreover, it is %known that under the former assumption, \eqref{eq: LD regular R^2} always holds if $B$ is either convex and open or is a %closed half-plane.
%\end{remark}

\subsection{Main proofs} \label{Sec: main proofs}
In this section we prove our main results Theorems~\ref{thm: perimeter LD},~\ref{thm: LD area}, and~\ref{thm: LD Gaussian area}. The proofs follow the same idea of using classical geometric inequalities to solve the variational problems~\eqref{eq: rate func I_P} and~\eqref{eq: rate func I_A} and thus find the rate functions $\mathcal J_P$ and $\mathcal J_A$.

\begin{proof}[\bf Proof of Theorem~\ref{thm: perimeter LD}.]
1.\ With a slight abuse of notation, denote by $P(h):= P(\conv(\im h))$ the perimeter of the convex hull of the image of a curve $h \in C_0[0,1]$.
%, i.e.\ the one-dimensional Lebesgue--Hausdorff measure of the boundary of the convex hull for the image of $h$. 
This is a continuous functional on $C_0[0,1]$, as follows from Cauchy's formula \eqref{eq: Cauchy}. By
\begin{equation} \label{eq: conve hull}
\conv (\{S_n(t) \}_{0 \le t \le 1})=\conv(S_0,S_1, \ldots, S_n),
\end{equation}
one has
\[
\frac12 P(S_n(\cdot)/n)= P_n/(2n).
\]
This equality, Mogulskii's LDP for trajectories of random walks (see Section~\ref{Sec: basics on LDP}), and the contraction principle (\cite[Theorem~4.2.1]{DemboZeitouni}) for continuous mappings yield that the sequence $(P_n/(2n))_{n \ge 1}$ satisfies an LDP in $\R$ with speed $n$ and the tight rate function
\begin{equation} \label{eq: rate function I_P}
\mathcal J_P(x):=\inf_{h \in C_0[0,1]: \, P(h) = 2x } I_C(h)=\min_{h \in AC_0[0,1]: \, P(h) = 2x } I_C(h).
\end{equation}
where, recall, $I_C$ is given by \eqref{eq: I def}. This implies \eqref{eq: rate func I_P}. We used that the lower semi-continuous non-negative function $I_C$ on $C_0[0,1]$ has compact sub-level sets and therefore it always attains its infimum over the closed set $\{P(h)=2x\}$.

The function $\mathcal J_P$ is lower semi-continuous on $\R$ as a rate function. It clearly satisfies $\mathcal J_P \leqslant \ubar I $. To show that it is strictly increasing on $[|\mu|, r_{max}]$, take any $x\neq |\mu|$ from this interval and choose an $h \in AC_0[0,1]$ such that $P(h)=2x$ and $\mathcal J_P(x)=I_C(h)$. If $h'(t) =h(1)$ a.e.\ $t$, then $\mathcal J_P(y) \le \ubar I(y) < \ubar I(x) = \mathcal J_P(x)$ for any $y \in [|\mu|, x)$, as required. Otherwise, for any $\varepsilon \in (0,1)$, 
\[
I_C((1-\varepsilon) h + \varepsilon h(1) \, \cdot) = \int_0^1 I((1-\varepsilon) h'(t) + \varepsilon h(1)) dt < \int_0^1 \big [(1-\varepsilon) I(h'(t) )+ \varepsilon I(h(1)) \big] dt,
\]
by strict convexity of $I$ on $\mathcal D_I$.  Hence $I_C((1-\varepsilon) h + \varepsilon h(1) \, \cdot) < I_C(h) = \mathcal J_P(x)$ by Jensen's inequality~\eqref{eqn: Jensen direct}. On the other hand, we have $P((1-\varepsilon) h + \varepsilon h(1) \,\cdot) <P(h) =2x$, which follows from Cauchy's formula~\eqref{eq: Cauchy} and the relation 
\[
\conv \big(\{(1-\varepsilon) h(t) + \varepsilon h(1) t\}_{0 \le t \le 1} \big) \subsetneq \conv(\{h(t)\}_{0 \le t \le 1}).
\]
The strict inequalities above imply strict monotonicity of $\mathcal J_P$ on $[|\mu|, r_{max}]$.

2.\ The main task is to find the minimum in \eqref{eq: rate function I_P} and its minimizers.

First consider the case $[0, |\mu|]$. 

For any function $h \in AC_0[0,1]$, it clearly holds $P(h) \ge 2|h(1)|$. Then by Jensen's inequality \eqref{eqn: Jensen direct} and the fact that $\ubar I$ is decreasing on $[0, |\mu|]$, for any $x \in [0, |\mu|]$. 
\[
\mathcal J_P (x)= \min_{h: P(h) = 2x} I_C(h) \ge \min_{h: P(h) = 2x} I(h(1)) \ge \min_{h: |h(1)| \le x} I(h(1)) = \min_{r \le x} \ubar{I}(r) = \ubar I(x).
\]
These inequalities are actually equalities since 
\begin{equation}
\label{eq: I_C(xlt)}
I_C( x \ell t) =   \ubar I (x), \quad P(x \ell t) = 2x, \qquad x >0, \ell \in \ubar{\Lambda}_x.
\end{equation}
Hence $J_P= \ubar I$ on $[0,|\mu|]$. Moreover, recalling that
$H_P(x)=\{h : I_C(h) = \ubar I(x), P(h) = 2x\}$,
\begin{equation} \label{eq: inf attained <}
%\bigl( I_C(h) = \ubar I(x), P(h) = 2x \bigr) \quad \text{ iff } \quad  h = x \ell_x %t, \qquad x \in (r_{min}, |\mu|].
H_P(x)=\{t \mapsto x \ell_x t\}, \qquad x \in (r_{min}, |\mu|].
\end{equation} 
Indeed, the facts that $\ubar I$ is strictly decreasing on $(r_{min}, |\mu|]$ and that Jensen's inequality~\eqref{eqn: Jensen direct} for the strictly convex rate function $I$ turns into equality only on functions with a.e.\ constant derivative, imply that the minimum in \eqref{eq: rate function I_P} is attained only on functions $h \in AC_0[0,1]$ that satisfy $P(h)= 2 |h(1)|=2x$, that is $h (t)= x \ell t$ for some $\ell \in \S^1$. The unique function $h$ of this form that satisfies the equality $I_C(h) = \ubar I(x)$ corresponds to the direction $\ell_x$.

The equality in \eqref{eq: LD per < E} for $x \in (r_{min}, |\mu|]$ now follows from the LDP for the perimeters $(P_n/(2n))_{n \ge 1}$ proved in Part 1. In fact, we have $\mathcal J_P  = \ubar I$ on $(r_{min}, |\mu|]$. On this interval $\ubar I$ is decreasing and convex (see Lemma~\ref{lem: properties of I_}.\ref{item: I_ unimodal}), hence continuous, and so the set $[0, x]$, which corresponds to the event $\{P_n \le 2xn\}$, is regular for the rate function $\mathcal J_P$. 

The claim in \eqref{eq: LD per < E} for $x = r_{min}$ holds trivially by $P_n \ge 2 r_{min} n$ a.s.

Consider now the case $[|\mu|, \infty)$.

Our main estimate follows from the inequality $I(v) \ge \conv \ubar I(|v|)$, $v \in \R^2$, and Jensen's inequality applied with the convex function $\conv \ubar I$. For any $h \in AC_0[0,1]$, we have
\begin{equation} 
\label{eq: I_C >=}
I_C(h) = \int_0^1 I(h'(t)) dt \ge \int_0^1 \conv \ubar I(|h'(t)|) dt \ge \conv \ubar I \Bigl( \int_0^1 |h'(t) | dt\Bigr) = \conv \ubar I \bigl( \Var(h) \bigr),
\end{equation}
where $\Var(h)$ denotes the total variation, i.e.\ the length, of a curve $h \in C_0[0,1]$.
%Note that $V(h)$ does not depend on the parametrization of the curve but $I_C(h)$ does.

Now use the following well-known inequality (see Corollary~\ref{cor:length of curve} in the Appendix), which is even referred to as geometric ``folklore'':  $\Var(h) \ge \frac12 P(h)$ for any $h \in C_0[0,1]$ of bounded variation. Since the function $\ubar I$  increases on $[|\mu|, r_{max}]$, so does its largest convex minorant $\conv \ubar I$. With the above, from \eqref{eq: I_C >=} we get: for $x \ge |\mu|$,
\begin{align} \label{eq: I_C >= 2}
\mathcal J_P(x)&=\min_{h: P(h) = 2x} I_C(h) \ge \min_{h: \Var(h) \ge x} I_C(h)  \notag \\
&\ge \min_{h: \Var(h) \ge x} \Bigl (\conv \ubar I \bigl( \Var(h) \bigr) \Bigr) \ge \min_{r \ge x} \bigl (\conv \ubar I \bigl( r \bigr) \bigr) = \conv \ubar I(x).
\end{align}
Using \eqref{eq: I_C(xlt)} for an upper bound, this gives $\conv \ubar I \leqslant \mathcal J_P \leqslant \ubar I $ on $[|\mu|, \infty)$.

We claim that if $\ubar I(x) = \conv \ubar I(x)$ for $x \ge |\mu|$, then
\begin{equation} \label{eq: inf attained >}
%\bigl( I_C(h) = \ubar I(x), P(h) = 2x \bigr ) \quad \text{ iff } \quad h \in \{x %\ell t, \ell \in \ubar{\Lambda}_x \}.
H_P(x)= \{t \mapsto x \ell t, \ell \in \ubar{\Lambda}_x \}.
\end{equation} 
We first note that by $\conv \ubar I \leqslant \mathcal J_P \leqslant \ubar I $  and the assumption $\ubar I(x) = \conv \ubar I(x)$, all inequalities in \eqref{eq: I_C >= 2} are  equalities. Then, since $\conv \ubar I$ is strictly increasing on $[|\mu|, r_{max})$, the infima in \eqref{eq: I_C >= 2} are attained  on the functions $h \in C_0[0,1]$ that satisfy $\Var(h) = \frac12 P(h)=x$. By Corollary~\ref{cor:length of curve} in the Appendix, such functions have the form $h(t) = |h(t)| \ell$ a.e.\ $t$ for some $\ell \in \S^1$ and satisfy $\Var(h)=x$. Further, the second inequality in \eqref{eq: I_C >=} is an equality iff $|h'(t)| \in [x_1, x_2]$ a.e.\ $t$, where $[x_1, x_2]$ is the maximal by inclusion interval that contains $x$ and is such that the restriction of $\conv \ubar{I}$ on $[x_1, x_2]$ is affine. Finally, the first inequality in \eqref{eq: I_C >=} is an equality for a function $h \in AC_0[0,1]$ that satisfies the conditions above iff 
\[
|h'(t)| \in \bigl\{y \in [x_1, x_2]: I(y \ell) = \conv \ubar{I}(y) \bigr\} =: L_x\text{ a.e. }t
\]
with the direction $\ell$ which was already fixed above. Since the rate function $I$ is strictly convex, so is $I(\cdot \ell)$, hence $L_x=\{x\}$. Thus we obtained that $|h'(t)| = x$ a.e.\ $t$ and by $I(x \ell) = \ubar I(x)$, we have $\ell \in \ubar{\Lambda}_x$. This finishes the proof of \eqref{eq: inf attained >}.

It remains to prove  \eqref{eq: LD per > E}.
In general, for an $x \in [ |\mu|, r_{max}]$  we can not assure regularity of the set $[x, \infty)$ (corresponding to the event $\{P_n \ge 2xn\}$) for the rate function $\mathcal J_P$.  The upper bound in~\eqref{eq: LD per > E} immediately follows from the LDP for the perimeters $(P_n/(2n))_{n \ge 1}$ we proved in Part 1 and the inequality $\conv \ubar I \leqslant \mathcal J_P$ (cf.~the upper bound in \eqref{eq: LD general} and~\eqref{eq: I_C >= 2}, respectively). For the lower bound in~\eqref{eq: LD per > E}, we consider two cases. If $\ubar I$ is continuous at $x$, then we use the inequality  $\mathcal J_P \leqslant \ubar I $ and the LDP for the perimeters (cf.~ the lower bound in~\eqref{eq: LD general}). If $\ubar I$ is discontinuous at $x$, then by Lemma~\ref{lem: properties of I_}.\ref{item: I_ discontinuous}, the distribution of $X_1$ has atoms at the points of $x \ubar{\Lambda}_x$, which must have equal weights  satisfying $\ubar I(x)= - \log \P(X_1 = x \ell)$ for $\ell \in \ubar{\Lambda}_x$. Then 
\[
\P(P_n \ge 2 x n) \ge \P(S_k = k x \ell, k=1, \ldots, n \,  \text{ for some } \ell \in \ubar{\Lambda}_x  )  = \#(\ubar{\Lambda}_x) e^{-n \ubar{I}(x)},
\]
which gives the lower bound in~\eqref{eq: LD per > E}.  The proof of \eqref{eq: LD per > E} is now finished.

3.\ The claims follow from the general statement \eqref{eq: LDP conditional} combined with  \eqref{eq: inf attained <}, \eqref{eq: inf attained >} and using that $\max_{0 \le k \le n} | S_k/n - h(k/n) \bigr| \le \max_{0 \le t \le 1} |S_n(t) - h(t)|$ for any $h \in C_0[0,1]$.
\end{proof}

\begin{proof}[\bf Proof of Theorem~\ref{thm: LD area}.]
Our argument is fully based on the ideas we developed in the proof of Theorem~\ref{thm: perimeter LD}. 

1. Denote by $A(h)$ the area of the convex hull of a curve $h \in C_0[0,1]$, i.e.\ $A(h):= A(\conv(\im h))$. It follows from the Steiner formula~\eqref{eq: Steiner}  that $A$ is a continuous functional on $C_0[0,1]$. From \eqref{eq: conve hull}, one has
\[
A(S_n(\cdot)/n)= A_n/n^2.
\]
This equality, Mogulskii's LDP for trajectories of random walks (see Section~\ref{Sec: basics on LDP}), and the contraction principle (\cite[Theorem~4.2.1]{DemboZeitouni}) for continuous mappings yield that the sequence $(A_n/n^2)_{n \ge 1}$ satisfies an LDP in $\R$ with speed $n$ and the tight rate function
\begin{equation} \label{eq: rate function I_A}
\mathcal J_A(a)=\inf_{h \in C_0[0,1]: \, A(h) = a } I_C(h) = \min_{h \in AC_0[0,1]: \, A(h) = a } I_C(h),
\end{equation}
where, recall, $I_C$ is given by \eqref{eq: I def}. This implies \eqref{eq: rate func I_A}. We used that the lower semi-continuous non-negative function $I_C$ on $C_0[0,1]$ has compact sub-level sets and therefore it always attains its infimum over the closed set $\{A(h)=a\}$. 

Clearly, $\mathcal J_A(a)=0$.  Let us check that $\mathcal J_A$ is strictly increasing on the set $\mathcal D_{\mathcal J_A}$. This assertion is trivial if this set is $\{0\}$, otherwise  pick a positive $a \in \mathcal D_{\mathcal J_A}$. Then $\mathcal J_A(a)=I_C(h)$ for some $h \in C_0[0,1]$ such that $A(h)=a$. Clearly, $h'$ is not constant a.e.\ on $[0,1]$ since otherwise $A(h)=0$. Consider the function $h_s$ such that $h_s(t)=(t/s)h(s) $ for $t \in [0,s]$ and $h_s=h$ on $[s,1]$, where $s \in (0,1]$; put $h_0:=h$. The area $A(h_s)$ decreases in $s$ and satisfies $A(h_0)=A(h)$, $A(h_1)=0$. By Jensen's inequality, we also have $I_C(h_s) \le I_C(h)$. From strict convexity of $I$, this inequality is strict if $A(h_s) < A(h)$. Since $A(h_s)$ is continuous in $s \in [0,1]$, these inequalities imply that $\mathcal J_A(a_1) < A(h) =\mathcal J_A(a)$ for any $a_1 \in [0, a)$.

Thus, the rate function $\mathcal J_A$ is left-continuous on $\mathcal D_{\mathcal J_A}$ since it is lower semi-continuous and increasing. Then
 \eqref{eq: LDP A} follows from (and is easily seen to be equivalent to) the LDP for the areas $(A_n/n^2)_{n \ge 1}$.
Finally, \eqref{eq: A shape} holds by the general result \eqref{eq: LDP conditional}.

2. The isoperimetric inequality for convex hulls,
\begin{equation} \label{eq: isoperimetric}
A(h) \le \Var(h)^2 /(2 \pi),
\end{equation}
is valid for any function $h \in C_0[0,1]$ of bounded variation. This is Ulam's version of the classical Dido problem, solved by Moran~\cite{Moran}. We have $\ubar I = \conv \ubar I$ by convexity of  $\ubar I$, which follows from rotational invariance of the distribution of $X_1$. Then by~\eqref{eq: I_C >=}, \eqref{eq: rate function I_A}, and \eqref{eq: isoperimetric}, 
\[
\mathcal J_A(a)= \min_{h \in C_0[0,1]: A(h) = a} I_C(h) \ge \min_{h: \Var(h)^2 \ge 2\pi a} I_C(h) \ge \min_{h: \Var(h) \ge \sqrt{2 \pi a}} \ubar I (\Var(h)) \ge \ubar I (\sqrt{2 \pi a})
\]
for $ a \ge 0$. These three inequalities actually are equalities, with the minima attained only on the functions that parametrize half circles with the constant speed $\sqrt{2 \pi a }$, and thus \eqref{eq: A shape invar} holds true. In fact, the value of $I_C$ on such a function $h$ is exactly $\ubar I (\sqrt{2 \pi a})$. Since $\ubar I (\sqrt{2 \pi a})$ strictly increases for $a \in [0, r_{max}]$, it must be that $a=A(h) = \Var(h)^2 /(2 \pi)$ and $|h'(t)|$ is constant for a.e.~$t$, ensuring that the second inequality in~\eqref{eq: I_C >=} is an equality. And the isoperimetric inequality \eqref{eq: isoperimetric}
is an equality only on parametrizations of semi-circles; see Tilli~\cite{Tilli}.

Finally, by Lemma~\ref{lem: properties of I_}.\ref{item: I_ unimodal}, the rate function $\mathcal J_A(a)$ is continuous on $[0, \infty)$, hence  \eqref{eq: LDP A} is valid for every $a \ge 0$. 
\end{proof}

\begin{proof}[\bf Proof of Theorem~\ref{thm: LD Gaussian area}.]
We need to find the rate function $\mathcal J_A$ given by \eqref{eq: rate func I_A}. For any $h \in AC_0[0,1]$, by Jensen's inequality we have
\[
I_C(h) = \frac12 \int_0^1 |h'(t)-\mu|^2 dt = \frac12 \int_0^1 |h'(t)|^2 dt  - h(1) \cdot \mu + \frac12 |\mu|^2 \ge \frac12 \Var(h)^2  - h(1) \cdot \mu + \frac12 |\mu|^2,
\]
where the inequality is an equality iff $|h'(t)|=\Var(h)$ for a.e.\ $t$. Hence, using that $\Var(h)$ is invariant under rotations of the image of $h$ about $0$, we get
\begin{equation} \label{eq: inf Gaussian shifted}
\mathcal J_A(a) = \min_{h \in C_0[0,1]: A(h) = a} I_C(h) = \min_{r \ge 0} \left( - r |\mu| + \min_{\substack{h: A(h) = a,\\ h(1)= r \mu / |\mu|}} \frac12 \Var(h)^2 \right)  + \frac12 |\mu|^2, \qquad a \ge 0.
\end{equation}

Assume $a>0$. It follows immediately from an approximation argument and the result by Pach~\cite{pach1978isoperimetric} for polygonal lines (see his Theorem~2 and the Remark just after it) that the above minimum over $h$ with the fixed endpoint $h(1)$ is attained only on parametrizations $h$ of circular arcs with $A(h)=a$. Denote by $R$ the radius of such an  arc and by $2 \varphi$ its angle, where $R>0$ and $0 \le \varphi \le \pi$. Then $\Var(h)= 2 \varphi R$, $\sin \varphi = r/(2 R)$, and $A(h)= \varphi R^2 - \frac12 r R \cos \varphi$ in both cases $0 \le \varphi \le \pi/2$ and $\pi/2 \le \varphi \le \pi$. Due to the fact that $\varphi / \sin \varphi$ is strictly increasing on $[0, \pi]$, the mapping $(\varphi, R) \mapsto (r, V)$ is a bijection between the sets $[0, \pi) \times (0,\infty)$ and $\{(r, V) \in [0, \infty) \times (0, \infty): r \le V \}$. Hence \eqref{eq: inf Gaussian shifted} reduces to
\begin{equation} \label{eq: minimization}
\mathcal J_A(a) =\frac12 |\mu|^2 + 2 \min_{\substack{0 \le \varphi \le \pi, R \ge 0: \\  R^2 (\varphi - \frac12 \sin 2 \varphi) = a}} \Bigl(  \varphi^2 R^2 - |\mu| R \sin \varphi \Bigr).
\end{equation}

Note that $R(\varphi)=\sqrt{a/(\varphi - \frac12 \sin 2 \varphi)}$ satisfies $\frac{\partial R}{\partial \varphi} = -a^{-1} R^3 \sin^2 \varphi$. The values of the function $\varphi^2 R(\varphi)^2 - |\mu| R(\varphi) \sin \varphi$ at $0$, $\pi/$, $\pi$ are respectively $+\infty$, $\pi a/2$, $\pi a$, hence this function attains its minimum at a critical point inside $(0, \pi)$ satisfying
\[
2 \varphi R^2 + 2 \varphi^2 R \frac{\partial R}{\partial \varphi} = |\mu| \frac{\partial R}{\partial \varphi} \sin \varphi + |\mu| R \cos \varphi.
\]
Dividing by $R^4$ and substituting the  expression for $\frac{\partial R}{\partial \varphi}$ gives 
\[
\frac{2 \varphi}{a}(\varphi - \sin \varphi \cos \varphi - \varphi \sin^2 \varphi) = \frac{|\mu|}{a R} (- \sin^3 \varphi + (\varphi - \sin \varphi \cos \varphi )\cos \varphi).
\]
Then $2 \varphi(\varphi  \cos^2 \varphi - \sin \varphi \cos \varphi) = \frac{|\mu|}{R} (- \sin \varphi + \varphi \cos \varphi )$,
and using that $\varphi \neq \tan \varphi $ on $(0, \pi/2)$,
\begin{equation} \label{eq: radius}
R = \frac{|\mu|}{2 \varphi \cos \varphi},
\end{equation}
which is possible only when $\varphi \in (0, \pi/2)$. This gives
\begin{equation}
\label{eq: a=}
\frac{a}{|\mu|^2} =  \frac{2\varphi - \sin 2 \varphi}{8 \varphi^2 \cos^2 \varphi}.
\end{equation}

It easy to check that this equation has only one solution $\varphi \in [0, \pi/2)$ for every $a\ge 0$. In fact, the right-hand side of \eqref{eq: a=} equals zero at $\varphi =0$ and $+\infty$ at $\varphi = \pi/2$, and its derivative
\[
\frac{1}{2\varphi^2 \cos^3\varphi} \left(
\cos^2 \varphi \sin \varphi + \varphi^2\sin \varphi - \varphi \cos \varphi
\right)
\]
is positive on $ (0, \pi/2)$ by
\begin{multline*}
	\cos^2 \varphi \sin \varphi + \varphi^2\sin \varphi - \varphi \cos \varphi>
	\cos^2 \varphi \sin \varphi + \sin^3 \varphi - \varphi \cos \varphi=\\
	\sin \varphi -\varphi \cos \varphi=
	\cos \varphi (\tan \varphi - \varphi)>0.
\end{multline*}

Substituting \eqref{eq: radius} into \eqref{eq: minimization} and using \eqref{eq: a=}, we obtain
\[
\mathcal J_A(a)  = \inf_{h \in C_0[0,1]: A(h) = a} I_C(h) = \frac{|\mu|^2}{2} \Bigl ( \frac{1}{\cos^2 \varphi} - \frac{2 \tan \varphi}{\varphi} + 1 \Bigr) = 4 \varphi  a -\frac12 |\mu|^2 \tan^2 \varphi.
\]
Clearly, this function is continuous in $a$, hence \eqref{eq: LDP A} is valid for every $a \ge 0$. The infimum is attained only at either of the two $\mu$-axially symmetric curves
\[
R \bigl( \sin( 2 \varphi t - \varphi) + \sin \varphi) ,\pm \cos( 2 \varphi t - \varphi) \mp \cos \varphi \bigr),
\]
where $\varphi$ and $R$ are given by \eqref{eq: radius} and \eqref{eq: a=}. This yields \eqref{eq: A shape shifted}.
\end{proof}

\subsection{The LDP's in continuous time} \label{sec: continuous time proof}
Here we obtain LDP's for convex hulls of L\'evy processes by reduction to random walks.

\begin{proof}[{\bf Proof of Theorem~\ref{thm: LD Levy}}]
First consider the perimeter $\mathsf P_T$ of the convex hull $\mathsf C_T= \conv(\{S_t\}_{0 \le t \le T}) $ of the L\'evy process $(S_t)_{t \ge 0}$. We shall compare it with the perimeter $P_{[T]}$ of the convex hull $C_{[T]}=\conv(0, S_1, \ldots, S_{[T]})$ of the random walk $(S_t)_{t \in \N \cup \{0\}}$. 

It follows from  Cauchy's formula~\eqref{eq: Cauchy} that
\[
0 \le (\mathsf P_T - P_{[T]})/(2 \pi) \le  \max_{k \in \{0, 1, \ldots, [T]\}} \sup_{k \le t \le k+1} |S_t - S_k| =:d_T, 
\]
where $d_T$ is an upper bound for the Hausdorff distance between $\mathsf C_T$ and $C_{[T]}$. Let us estimate probabilities of large deviations of $d_T/T$. By stationarity of increments of $(S_t)_{t \ge 0}$, for every $\varepsilon>0$ we have
\[
\P \big (d_T\ge \varepsilon T \big ) \le ([T] +1) \P \Big( \sup_{0 \le t \le 1} |S_t | \ge \varepsilon T  \Big). 
\]

Put $\tilde S_t :=S_t - t \mu$ for $t \ge 0$ (recall that $S_1 = X_1$) and let $\tilde S_t^{(1)}$ and $\tilde S_t^{(2)}$ be the coordinates of $\tilde S_t$ is any orthonormal basis of $\R^2$. Note that
\[
\sup_{0 \le t \le 1} |S_t | \le |\mu| + \sup_{0 \le t \le 1} |\tilde S_t | \le |\mu| + \max_{i, j \in \{1, 2\} } \sup_{0 \le t \le 1} \big( (-1)^i \tilde S_t^{(j)} \big).
\]
Denote $a:= \varepsilon T- |\mu|$. Then for any $u>0$ we get
\[
\P \Big( \sup_{0 \le t \le 1} |S_t | \ge \varepsilon T   \Big)  \le  \sum_{i, j =1}^2  \P \Big(  \sup_{0 \le t \le 1} \big( (-1)^i \tilde S_t^{(j)} \big) \ge a \Big)  = \sum_{i, j =1}^2 \P \Big(  \sup_{0 \le t \le 1} e^{u (-1)^i \tilde S_t^{(j)}} \ge e^{u a} \Big). 
\]

Since $(\tilde S_t)_{t \ge 0}$ is a zero-mean L\'evy process in $\R^2$, each of the four stochastic processes $( (-1)^i \tilde S_t^{(j)})_{t \ge 0}$ is a right-continuous real-valued martingale. Then $ (e^{u (-1)^i \tilde S_t^{(j)}})_{t \ge 0}$ are right-continuous positive sub-martingales, because $x \mapsto e^{u x}$ is a positive convex function of $ x \in \R$. Hence, applying Doob's maximal inequality (Revuz and Yor~\cite[Chapter II, Theorem~1.7]{RevuzYor}), we obtain
\[
\P \Big( \sup_{0 \le t \le 1} |S_t | \ge \varepsilon T   \Big)  \le \sum_{i, j =1}^2  e^{-ua } \E e^{u (-1)^i \tilde S_1^{(j)}} = \sum_{i, j =1}^2  \exp \Big\{-  \Big( ua - \log \E e^{u (-1)^i \tilde S_1^{(j)}}  \Big) \Big\}. 
\]
Finally, if $a>0$ (where $0=\E \tilde S_1^{(j)} $), then optimizing the last expression over $u >0$ yields
\[
\P \big (d_T \ge \varepsilon T \big ) \le ([T] +1)  \sum_{i, j =1}^2 \exp \big \{- I_{i, j} \big (\varepsilon T - |\mu| \big) \big\},
\]
where $I_{i,j}$ denotes the rate function of $(-1)^i \tilde S_1^{(j)}$.

Since the Laplace transform of $S_1$ is finite in $\R^2$ by the assumption, the Laplace transform of each of the random variables $(-1)^i \tilde S_1^{(j)}$ is finite in $\R$. This implies $\lim_{u \to \infty} I_{i,j}(u)/u = \infty$; see~Rockafellar~\cite[Theorems~8.5 and 13.3]{Rockafellar} or Vysotsky~\cite[Eqs.~(5.4) and (5.5)]{VysotskyNote}. Therefore, for every $\varepsilon >0$, we have
\begin{equation} \label{eq: d_T exp zero}
\lim_{T \to \infty} \frac{1}{T} \log \P \big (d_T/T> \varepsilon \big) = -\infty,
\end{equation}
which means that the sequence of random variables $(d_T/T)_{T >0} $ is {\it exponentially equivalent} to $0$ as $T \to \infty$ in the sense of Definition~4.2.10 in~\cite{DemboZeitouni}. 

Finally, let us use that
\[
|\mathsf P_T/T - P_{[T]}/[T] | \le |\mathsf P_T/T -  P_{[T]}/T | + P_{[T]}(1/[T]-1/T) \le 2 \pi d_T/T + (P_{[T]}/[T])/T, 
\]
where the r.h.s.\ is exponentially equivalent to $0$ as $T \to \infty$ by \eqref{eq: d_T exp zero} and the fact that $  (P_n/n)_{n \in \N} $  satisfies an LDP in $\R$ with a tight rate function (by Theorem~\ref{thm: perimeter LD}). Therefore,
the sequences $\mathsf P_T/(2T) $ and $  P_{[T]}/(2[T]) $ are  exponentially equivalent as $T \to \infty$, hence they satisfy the same LDP by \cite[Theorem~4.2.13]{DemboZeitouni}, as claimed.

As for the areas, the Steiner formula~\eqref{eq: Steiner} yields
\[
0 \le (\mathsf A_T - A_{[T]})  \le  P_{[T] } d_T + \pi d_T^2,
\]
and it follows by the same argument as above that $\mathsf A_T/T^2 $ and $  A_{[T]}/T^2 $ are exponentially equivalent  as $T \to \infty$ (use \eqref{eq: d_T exp zero} and the facts that $  (P_n/n)_{n \in \N} $ and $  (A_n/n^2)_{n \in \N} $ satisfy LDPs  with tight rate functions). Then it follows from Theorem~\ref{thm: LD area} and \cite[Theorem~4.2.13]{DemboZeitouni} that $(\mathsf A_T/T^2)_{T\ge1} $ and $ ( A_{[T]}/[T]^2)_{T \ge1} $ satisfy the same LDP, as claimed.
\end{proof}

\subsection{The LDP's under the Cram\'er moment assumption} \label{sec: LDP Cramer}
Here we partially extend our main Theorems~\ref{thm: perimeter LD} and~\ref{thm: LD area} under the weaker assumption $0 \in \intr \mathcal D_{\mathcal L}$. We will use the contraction principle by~Vysotsky~\cite{VysotskyNote}.

Denote by $BV[0,1]=BV([0,1];\R^2)$ the set of right-continuous functions of bounded variation from $[0,1]$ to $\R^2$. Denote by
$A(h)$ and $P(h)$  respectively the area and the perimeter of $\conv(h([0,1]) \cup \{ 0 \}) $ of an $h \in BV[0,1]$. This extends the definitions given in Section~\ref{Sec: main proofs} for $h \in C_0[0,1]$. Consider the functional
\begin{equation} \label{eq: I BV} 
I_{BV}(h):= \sup_{\mathbf t \subset (0,1]: \, \# \mathbf t < \infty } I_C(h^{\mathbf t} ), \quad h \in BV[0,1],
\end{equation}
where $h^{\mathbf t}$ denotes the continuous function on $[0,1]$ defined by linear interpolation between its values at $\mathbf t \cup \{0,1\}$ that are given by $h^{\mathbf t}(s) := h(s)$ for $s \in \mathbf t \cup \{1\}$ and $h^{\mathbf t}(0):=0$. This functional satisfies $I_{BV}=I_C$ on $AC_0[0,1]$; see~\cite[Theorem~5.1]{VysotskyNote}, which gives an explicit and transparent formula for $I_{BV}(h)$ in terms of the Lebesgue decomposition of $h$.

\begin{prop} \label{prop: Cramer}
Assume that $X_1$ is a random vector in the plane such that $0 \in \intr \mathcal D_{\mathcal L}$. Then the  random variables $(P_n/(2n))_{n \ge 1}$ and  $(A_n/n^2)_{n \ge 1}$ satisfy the LDP's in $\R$ with speed $n$ and  the respective tight rate functions $\tilde{\mathcal J}_P $ and $\tilde{\mathcal J}_A$ given by
\begin{equation} \label{eq: rate func Cramer}
\tilde{\mathcal J}_P(x):=\cl \inf_{\substack{h \in BV[0,1]: \\ P(h) = 2x }} I_{BV}(h), \qquad \tilde{\mathcal J}_A(x):=\cl \inf_{\substack{h \in BV[0,1]: \\ A(h) = x}} I_{BV}(h), \quad x \ge 0.
\end{equation}
These rate functions increase on $[|\mu|, \infty)$ and $[0, \infty)$, respectively. We always have $\tilde{\mathcal J}_P = \ubar I$ on $[0, |\mu|]$. Moreover, $\tilde{\mathcal J}_P = \ubar I$ if $\ubar I$ is convex. Also, we have $\tilde{\mathcal J}_A(a) = \ubar I (\sqrt{2 \pi x})$ for $x \ge 0$ if the distribution of $X_1$ is rotationally invariant.
\end{prop}

Note that the monotonicity properties of $\tilde{\mathcal J}_P$ and $\tilde{\mathcal J}_A$ imply that the lower semi-continuous regularizations $\cl$ in \eqref{eq: rate func Cramer} may change the values of the infima only at the discontinuity points.

\begin{proof}
Let us equip $BV[0,1]$ with the metric $\rho$ equal the Hausdorff distance between the {\it completed graphs} of functions, defined by  $\Gamma h:= \{(t, x): 0 \le t \le 1, x \in  [h(t-), h(t)] \}$ for $h \in BV[0,1]$, where $h(0-):=0$. Note that $\Gamma h$ is a compact subset of $[0,1] \times \R^2$ and it uniquely defines $h$, i.e.\ $\Gamma h_1=\Gamma h_2$ for $h_1, h_2 \in BV[0,1]$ implies $h_1=h_2$. The total variation of an $h \in BV[0,1]$, given by $\Var(h):= \sup_{\mathbf t \subset (0,1]: \, \# \mathbf t < \infty } \Var(h^{\mathbf t} )$, is simply the length of the spatial coordinate of any continuous bijective parametrization of~$\Gamma h$.

It follows from Steiner's and Cauchy's formulas \eqref{eq: Steiner} and \eqref{eq: Cauchy} that the functionals $A$ and $P$ are continuous in the metric $\rho$ and moreover, they are {\it uniformly} continuous on the sets $\{h \in BV[0,1]: \Var(h)\le R\}$ for every $R>0$. Therefore, the contraction principle for the trajectories $S_n(\cdot)$ in $BV[0,1]$, given by Theorem~3.3 in~\cite{VysotskyNote}  (which uses a metric longer than $\rho$, see~\cite[Eqs.~(2.6) and (2.7)]{VysotskyNote}),  yields the LDPs stated with the respective rate functions given in~\eqref{eq: rate func Cramer}.
%Let us comment that this contraction principle is needed to bring into a standard form the weaker LDP-type %result by Borovkov and Mogulskii~\cite{BorovkovMogulskii2} for trajectories of random walks under the Cram%\'er moment assumption.

The rest of the proof is identical to the ones of the corresponding parts of Theorems~\ref{thm: perimeter LD} and~\ref{thm: LD area}. We comment only on the differences. The monotonicity properties  of $\tilde{\mathcal J}_P$ and $\tilde{\mathcal J}_A$ follow from equalities \eqref{eq: I BV} and \eqref{eq: rate func Cramer}. We get only non-strict monotonicity since we are not claiming that the infima in \eqref{eq: rate func Cramer} are always attained, as opposed to the main case $\mathcal D_{\mathcal L}=\R^2$.

Furthermore, by \eqref{eq: I BV} and Jensen's inequality, we have $I_{BV}(h) \ge I(h(1))$ for any $h \in BV[0,1]$. Moreover, if $\ubar I$ is convex, we have $I_{BV}(h) \ge \ubar I(\Var(h))$. This follows from  \eqref{eq: I_C >=} and \eqref{eq: I BV} using lower semi-continuity of $\ubar I$ (Lemma~\ref{lem: properties of I_}.\ref{item: I_ unimodal}) if we choose an increasing sequence $({\mathbf t}_n)_{n \ge 1}$ of finite subsets of $[0,1]$ such that $I_C(h^{{\mathbf t}_n}) \to I_{BV}(h)$ and $\Var(h^{{\mathbf t}_n}) \to \Var(h)$ as $n \to \infty$. The  two inequalities above for $I_{BV}(h)$ yield, as in the proof of Theorem~\ref{thm: perimeter LD}, that $\inf_{h: P(h)=2x} I_{BV}(h) \ge \ubar I(x)$ for any $x \in [0, |\mu|]$ and also for $x \ge |\mu|$ if $\ubar I$ is convex. Hence $\tilde{\mathcal J}_P(x) \ge \ubar I(x)$ for such $x$ since  $\ubar I$ is lower semi-continuous.  On the other hand, we have
\[
\tilde{\mathcal J}_P(x) \le \inf_{\substack{h \in BV[0,1]: \\P(h) = 2x }} I_{BV}(h) \le \inf_{\substack{h \in AC_0[0,1]: \\ P(h) = 2x }} I_{BV}(h) = \inf_{\substack{h \in AC_0[0,1]: \\ P(h) = 2x }} I_C(h) \le \ubar I(x),
\]
where we used that $I_{BV} = I_C$ on $AC_0[0,1]$. This yields the claims on $\tilde{\mathcal J}_P $.

Similarly, if $\ubar I$ is convex, which is surely the case when  the distribution of $X_1$ is rotationally invariant, then we have $I_{BV}(h) \ge \ubar I(\Var(h))$ for $h \in BV[0,1]$, hence  $\inf_{h: A(h)=a} I_{BV}(h) \ge \ubar I (\sqrt{2 \pi a})$ for $a \ge 0$ by  the same argument as in the proof of Theorem~\ref{thm: LD area}. Hence $\tilde{\mathcal J}_A(a) \ge  \ubar I (\sqrt{2 \pi a})$ for $a \ge 0$. On the other hand,  for rotationally invariant distributions of $X_1$ we have $\tilde{\mathcal J}_A(a) \le  \ubar I (\sqrt{2 \pi a})$, arguing as above for $\tilde{\mathcal J}_P(x) \le \ubar I(x)$. This yields the claim on~$\tilde{\mathcal J}_A $.
\end{proof}

\section*{Appendix}
\subsection*{Perimeters} Throughout the paper, by the {\it perimeter} $P(C)$ of a non-empty convex set $C$ on the plane we mean the length of its boundary unless $C$ is a line segment, in which case $P(C)$ is its doubled length. 
Recall that a continuous curve in $\R^d$ is {\it rectifiable} if it has finite length (equivalently, it has bounded variation). 

The following simple proposition is proved %[[in a slightly weaker form]] 
in our separate note~\cite{AkopyanVysotskyGeometry}, which was initially motivated by the questions concerning the perimeter of the convex hulls considered in the present paper. For the reader's convenience, we present the result here. Its main use here is in the corollary, which not only gives the ``folklore'' inequality for the half-perimeter but also specifies all instances when the equality is attained.

\begin{prop}
   \label{prop:length of curve}
Let $\gamma$ be a rectifiable curve in $\R^2$, %[[ with the endpoints $a$ and $b$]] 
and let $\Gamma$ denote its convex hull. Then
\[
\length \gamma \geq \per \Gamma - \diam \Gamma. %|b-a|.
\]
\end{prop}	
\begin{cor}
   \label{cor:length of curve}
It holds that
\[
\length \gamma \geq \frac12 \per \Gamma,
\]
and equation can be attained only if  $\gamma$ parametrizes is a line segment.
\end{cor}
\begin{remark} \label{rem: mean width}
These statements remain valid if we replace $\R^2$ by $\R^d$ (with any $d \ge 2$) and $\per \Gamma$ by $\frac{d v_d}{v_{d-1}} W(\Gamma)$, where $v_d$ denotes volume of a unit ball in $\R^d$ and 
\begin{equation} \label{eq: mean width}
W(\Gamma):= \frac{1}{|\S^{d-1}|} \int_{\S^{d-1}} w_\ell(\Gamma) d \ell
\end{equation}
is {\it mean width} of $\Gamma$, with $w_\ell (\Gamma)$ being width of $\Gamma$ in the direction $\ell$ i.e.\ length of the projection of $\Gamma$ on the line passing through the origin in the direction $\ell$. The normalizing factor corresponds to mean width $\frac{2 v_{d-1}}{d v_d}$ of a unit segment in $\R^d$.
\end{remark}

It is easy to prove the remark using {\it Crofton's formula} (Schneider and Weil~\cite[Eq.~(5.32)]{SchneiderWeil})
\[
		 \length \gamma = \frac{1}{v_{d-1}} \iint \limits_{\mathbb{S}^{d-1} \mathbb{R_+}}n_{\gamma}(\ell, r) d\ell dr, 
\]
where $n_\gamma(\ell, r)$ denotes the number of intersections of $\gamma$ with the hyperplane perpendicular to the direction $\ell$ at the distance $r$ from the origin. Indeed, consider the closed curve $\gamma'$ obtained by joining the end points of $\gamma$ by a line segment. Almost every hyperplane intersecting $\Gamma$ intersect $\gamma'$ at least at two points since $\conv(\gamma')=\Gamma$. It remains to use that $|\S^{d-1}| =  d v_d$. 

Note that  Crofton's formula implies   {\it Cauchy's formula} for the perimeter of the planar convex set $\Gamma$:
	\begin{equation} \label{eq: Cauchy}
		 \per \Gamma =  \frac12 \int_{\S^1} w_\ell(\Gamma) d\ell.
	\end{equation}
	
\subsection*{Areas}
Let $C \subset \R^2$ be a non-empty bounded convex set and let $B \subset \R^2$ be the closed unit ball centred at the origin. {\it Steiner's  formula} (\cite[Eq.~(14.5)]{SchneiderWeil}) asserts that for every $r>0$,
\begin{equation} \label{eq: Steiner}
A(C+r B) = A(C) +  P( C) r + \pi r^2, 
\end{equation}	
where `$+$' stands for Minkowski addition of sets. 

\subsection*{Measurability}
Let us show that the perimeters and areas $(P_n)_{n \in \N}, (A_n)_{n \in \N}$, $(\mathsf P_T)_{T >0}$, $(\mathsf A_T)_{T >0}$ of the convex hulls, introduced in Sections~\ref{Sec: intro} and~\ref{sec: continuous time}, are measurable.

It follows from \eqref{eq: Cauchy} and \eqref{eq: Steiner} that for every $n \in \N$, the mappings 
\[
(x_1, \ldots, x_n) \mapsto P(\conv(0, x_1, \ldots, x_n )) \quad \text{and} \quad (x_1, \ldots, x_n) \mapsto A(\conv(0, x_1, \ldots, x_n ))
\]
are continuous from $\R^{2 \times n} $ to $\R$. Hence $P_n$ and $A_n$ are random variables.

Furthermore, for any $T>0$ and a dense subset $\{t_k\}_{k \in \N}$ of $[0, T]$ that includes $T$,
\begin{align*}
\cl \mathsf C_T &= \cl ( \conv(\{S_t\}_{0 \le t \le T}) )=  \conv ( \cl (\{S_t\}_{0 \le t \le T}) ) \\
&= \conv (  \cl (\{S_{t_k}\}_{k \in \N}) ) = \cl (  \conv (\{S_{t_k}\}_{k \in \N}) ) \quad \text{a.s.,}
\end{align*}
where  the second and the fourth equalities hold true by~\cite[Theorem~17.2]{Rockafellar}, which applies because the trajectories of a L\'evy process are bounded a.s.\ on any interval, and in the third equality we used that the trajectories are right-continuous and have left limits a.s. Then $\cl \mathsf C_T=   \cl (  \cup_{k=1}^\infty \conv ( S_{t_1} \ldots, S_{t_k}) )$  by Carath\'eodory's theorem (\cite[Theorem~17.1]{Rockafellar}).
 
Hence, since the union on the r.h.s.\ is a convex set, we have 
\[
\mathsf A_T=A(\mathsf C_T)= A (  \cup_{k=1}^\infty \conv ( S_{t_1} \ldots, S_{t_k}) ) = \lim_{k \to \infty}A( \conv ( S_{t_1} \ldots, S_{t_k})),
\]
and by the above, $\mathsf A_T$ is measurable as a limit of measurable functions. Also, for any  $\ell \in \S^1$,
\[
w_\ell (\mathsf C_T)= w_\ell ( \cup_{k=1}^\infty \conv ( S_{t_1} \ldots, S_{t_k}) )= \lim_{k \to \infty} w_\ell ( \conv ( S_{t_1} \ldots, S_{t_k})),
\]
which yields $\mathsf P_T = \lim_{k \to \infty}P( \conv ( S_{t_1} \ldots, S_{t_k}))$ by \eqref{eq: Cauchy} and the monotone convergence theorem. Hence $\mathsf P_T $ is measurable as a limit of measurable functions.

\section*{Acknowledgements}  We are grateful to Andrew Wade for bringing the perimeter problem to our attention, and to Endre Makai for referring us to the paper~\cite{pach1978isoperimetric} by J{\'a}nos Pach. We wish to thank Fedor Petrov for showing us a simple proof of Proposition~\ref{prop: generalization}.\ref{item: conv I_ =}. We are indebted to the anonymous referees for their comments and the suggestion to include continuous time results.

\bibliographystyle{plain}
\bibliography{convdev}

\begin{thebibliography}{10}

\bibitem{AkopyanVysotskyGeometry}
Arseniy Akopyan and Vladislav Vysotsky.
\newblock On the lengths of curves passing through boundary points of a planar
  convex shape.
\newblock {\em Amer. Math. Monthly}, 124:588--596, 2017.

\bibitem{AKMV}
Gerold Alsmeyer, Zakhar Kabluchko, Alexander Marynych, and Vladislav Vysotsky.
\newblock How long is the convex minorant of a one-dimensional random walk?
\newblock {\em Electron. J. Probab.}, 25:1--22, 2020.

\bibitem{B-N}
Ole Barndorff-Nielsen.
\newblock {\em Information and exponential families in statistical theory}.
\newblock John Wiley \& Sons, Ltd., Chichester, 1978.

\bibitem{BNielsen}
Ole Barndorff-Nielsen and Glen Baxter.
\newblock Combinatorial lemmas in higher dimensions.
\newblock {\em Trans. Amer. Math. Soc.}, 108:313--325, 1963.

\bibitem{Bertoin}
Jean Bertoin.
\newblock {\em L\'{e}vy processes}.
\newblock Cambridge University Press, Cambridge, 1996.

\bibitem{BorovkovMogulskii2}
Aleksandr~A. Borovkov and Anatolii~A. Mogulskii.
\newblock Large deviation principles for random walk trajectories. {II}.
\newblock {\em Theory Probab. Appl.}, 57:1--27, 2013.

\bibitem{Claussen+}
Gunnar Claussen, Alexander~K. Hartmann, and Satya~N. Majumdar.
\newblock Convex hulls of random walks: large-deviation properties.
\newblock {\em Phys. Rev. E.}, 91:052104, 2015.

\bibitem{Croft+}
Hallard~T. Croft, Kenneth~J. Falconer, and Richard~K. Guy.
\newblock {\em Unsolved problems in geometry}.
\newblock Springer-Verlag, New York, 1991.

\bibitem{DemboZeitouni}
Amir Dembo and Ofer Zeitouni.
\newblock {\em Large deviations techniques and applications}.
\newblock Springer-Verlag, Berlin, 2010.
\newblock Corrected reprint of the second (1998) edition.

\bibitem{glaeser2016universe}
Georg Glaeser, Hellmuth Stachel, and Boris Odehnal.
\newblock {\em The Universe of Conics: From the ancient Greeks to 21st century
  developments}.
\newblock Springer, 2016.

\bibitem{Khoshnevisan}
Davar Khoshnevisan.
\newblock Local asymptotic laws for the {B}rownian convex hull.
\newblock {\em Probab. Theory Related Fields}, 93:377--392, 1992.

\bibitem{KuelbsLedoux}
James Kuelbs and Michel Ledoux.
\newblock On convex limit sets and {B}rownian motion.
\newblock {\em J. Theoret. Probab.}, 11:461--492, 1998.

\bibitem{McRedmondWade}
James McRedmond and Andrew~R. Wade.
\newblock The convex hull of a planar random walk: perimeter, diameter, and
  shape.
\newblock {\em Electron. J. Probab.}, 23:Paper No. 131, 1--24, 2018.

\bibitem{Mogulskii}
Anatolii~A. Mogulskii.
\newblock Large deviations for the trajectories of multidimensional random
  walks.
\newblock {\em Theor. Probab. Appl.}, 21(2):300--315, 1976.

\bibitem{MolchanovWespi}
Ilya Molchanov and Florian Wespi.
\newblock Convex hulls of {L}\'{e}vy processes.
\newblock {\em Electron. Commun. Probab.}, 21:Paper No. 69, 11, 2016.

\bibitem{Moran}
P.~A.~P. Moran.
\newblock On a problem of {S}. {U}lam.
\newblock {\em J. London Math. Soc.}, 21:175--179, 1946.

\bibitem{pach1978isoperimetric}
J{\'a}nos Pach.
\newblock On an isoperimetric problem.
\newblock {\em Studia Sci. Math. Hungar.}, 13:43--45, 1978.

\bibitem{RevuzYor}
Daniel Revuz and Marc Yor.
\newblock {\em Continuous martingales and {B}rownian motion}.
\newblock Springer-Verlag, Berlin, third edition, 1999.

\bibitem{Rockafellar}
R.~Tyrrell Rockafellar.
\newblock {\em Convex analysis}.
\newblock Princeton University Press, Princeton, N.J., 1970.

\bibitem{SchneiderWeil}
Rolf Schneider and Wolfgang Weil.
\newblock {\em Stochastic and integral geometry}.
\newblock Springer-Verlag, Berlin, 2008.

\bibitem{SnyderSteele}
Timothy~Law Snyder and J.~Michael Steele.
\newblock Convex hulls of random walks.
\newblock {\em Proc. Amer. Math. Soc.}, 117:1165--1173, 1993.

\bibitem{SpitzerWidom}
Frank Spitzer and Harold Widom.
\newblock The circumference of a convex polygon.
\newblock {\em Proc. Amer. Math. Soc.}, 12:506--509, 1961.

\bibitem{Tilli}
Paolo Tilli.
\newblock Isoperimetric inequalities for convex hulls and related questions.
\newblock {\em Trans. Amer. Math. Soc.}, 362(9):4497--4509, 2010.

\bibitem{VysotskyStrict}
Vladislav Vysotsky.
\newblock When is the rate function of a random vector strictly convex?
\newblock {\em Preprint}, 2020.
\newblock Available at arXiv:2009.06809 [math.PR].

\bibitem{VysotskyNote}
Vladislav Vysotsky.
\newblock Contraction principle for trajectories of random walks and
  {C}ram\'er's theorem for kernel-weighted sums.
\newblock {\em Accepted in ALEA Lat. Am. J. Probab. Math. Stat.}, 2021.
\newblock Available at arXiv:1909.00374 [math.PR].

\bibitem{VysotskyZaporozhets}
Vladislav Vysotsky and Dmitry Zaporozhets.
\newblock Convex hulls of multidimensional random walks.
\newblock {\em Trans. Amer. Math. Soc.}, 370:7985--8012, 2018.

\bibitem{WadeXu}
Andrew~R. Wade and Chang Xu.
\newblock Convex hulls of planar random walks with drift.
\newblock {\em Proc. Amer. Math. Soc.}, 143:433--445, 2015.

\bibitem{WadeXu2}
Andrew~R. Wade and Chang Xu.
\newblock Convex hulls of random walks and their scaling limits.
\newblock {\em Stochastic Process. Appl.}, 125(11):4300--4320, 2015.

\end{thebibliography}

\end{document}